\documentclass[11pt]{article}
\usepackage{latexsym,amssymb}
\usepackage{amsmath}

\newtheorem{Theorem}{\bf Theorem}[section]
\newtheorem{Lemma}{\bf Lemma}[section]
\newtheorem{Proposition}{\bf Proposition}[section]
\newtheorem{Corollary}{\bf Corollary}[section]
\newtheorem{Remark}{\bf Remark}[section]
\newtheorem{Example}{\bf Example}[section]
\newtheorem{Definition}{\bf Definition}[section]

\newenvironment{theorem}{\begin{Theorem}$\!\!\!$}{\end{Theorem}}
\newenvironment{lemma}{\begin{Lemma}$\!\!\!$}{\end{Lemma}}
\newenvironment{proposition}{\begin{Proposition}$\!\!\!$}{\end{Proposition}}
\newenvironment{corollary}{\begin{Corollary}$\!\!\!$}{\end{Corollary}}
\newenvironment{remark}{\begin{Remark}$\!\!\!$}{\end{Remark}}

\numberwithin{equation}{section}

\begin{document}

\title{Asymptotics for a nonlinear integral equation\\ with a generalized heat kernel}
\author{Kazuhiro Ishige\\
Mathematical Institute, Tohoku University\\
Aoba, Sendai 980-8578, Japan,
\\
\\
Tatsuki Kawakami\\
Department of Mathematical Sciences\\
Osaka Prefecture University\\
Sakai 599-8531, Japan
\\
\\
and
\\
\quad
\\
Kanako Kobayashi\\
Mathematical Institute, Tohoku University\\
Aoba, Sendai 980-8578, Japan}
\date{}
\maketitle
\begin{abstract}
This paper is concerned with a nonlinear integral equation 
$$
(P)\qquad
u(x,t)=\int_{{\bf R}^N}G(x-y,t)\varphi(y)dy+\int_0^t\int_{{\bf R}^N}G(x-y,t-s)f(y,s:u)dyds, 
\quad
$$
where $N\ge 1$, $\varphi\in L^\infty({\bf R}^N)\cap L^1({\bf R}^N,(1+|x|^K)dx)$ for some $K\ge 0$. 
Here $G=G(x,t)$ is a generalization of the heat kernel. 
We are interested in the asymptotic expansions of the solution of $(P)$ 
behaving like a multiple of the integral kernel $G$ as $t\to\infty$. 
\end{abstract}
\vspace{10pt}
\newpage
\section{Introduction}
Let $u$ be a solution of a nonlinear integral equation, 
\begin{equation}
u(x,t)=\int_{{\bf R}^N}G(x-y,t)\varphi(y)dy
+\int_0^t\int_{{\bf R}^N}G(x-y,t-s)f(y,s:u)dyds
\label{eq:1.1}
\end{equation}
in ${\bf R}^N\times(0,\infty)$, where $N\ge 1$, 
$f$ is an inhomogeneous term possibly depending on the solution $u$ itself, 
and $G=G(x,t)$ is an integral kernel satisfying the following condition:  
\begin{itemize}
  \item[({\bf G})]
  \begin{itemize}
  	\item[(i)] $G\in C^\gamma({\bf R}^N\times(0,\infty))$ for some $\gamma\in{\bf N}$;
	\item[(ii)] 
	There exist positive constants $d$ and $L$ such that 
  	\begin{eqnarray}
	\label{eq:1.2}
	 & & G(x,t)=t^{-\frac{N}{d}}G\biggr(\frac{x}{t^{1/d}},1\biggr),
	\qquad x\in{\bf R}^N,\,\,\,t>0,\\
	\label{eq:1.3}
	 & & \sup_{x\in{\bf R}^N}\,(1+|x|)^{N+L+j}|\nabla_x^j G(x,1)|<\infty,
	 \qquad j=0,\dots,\gamma;\qquad\qquad\qquad
	 \end{eqnarray}
	\item[(iii)]\,\,
	$\displaystyle{
	G(x,t)=\int_{{\bf R}^N}G(x-y,t-s)G(y,s)dy}
	$\quad
	for $x\in{\bf R}^N$ and $t>s>0$.
  \end{itemize}
\end{itemize}
Condition~$({\bf G})$ holds for the fundamental solutions of the following linear diffusion equations, 
$$
\begin{array}{ll}
\partial_t u+(-\Delta)^{\theta/2}u=0
 & \mbox{in}\quad{\bf R}^N\times(0,\infty)\qquad (0<\theta<2),\vspace{5pt}\\
\partial_t u+(-\Delta)^mu=0
 & \mbox{in}\quad{\bf R}^N\times(0,\infty)\qquad (m=1,2,3,\dots),
\end{array}
$$
and integral equation~\eqref{eq:1.1} appears in the study of various nonlinear diffusion equations. 
In this paper we give the asymptotic expansions 
of the solutions of \eqref{eq:1.1} behaving like a multiple of the kernel $G$ as $t\to\infty$. 
Our arguments are applicable to the large class of nonlinear diffusion equations, including 
the following semilinear parabolic equations (see Section~6): 
\begin{itemize}
  \item (Fractional semilinear parabolic equation)
  \begin{equation}
  \label{eq:1.4}
  \partial_t u+(-\Delta)^{\theta/2}u=|u|^{p-1}u\quad\mbox{in}\quad{\bf R}^N\times(0,\infty),
  \end{equation}
  where $N\ge 1$, $0<\theta<2$ and $p>1+\theta/N$ 
  (see e.g. \cite{AF}, \cite{FIK}, \cite{FK}, \cite{HKN} and \cite{S});
  \item (Higher order semilinear parabolic equation)
  \begin{equation}
  \label{eq:1.5}
  \partial_t u+(-\Delta)^m u=|u|^p\quad\mbox{in}\quad{\bf R}^N\times(0,\infty),
  \end{equation} 
  where $N\ge 1$, $m=1,2,\dots$ and $p>1+2m/N$
  (see e.g. \cite{CM}, \cite{EGKP1}, \cite{EGKP2}, \cite{GP} and \cite{GG}).
\end{itemize}
See also a forthcoming paper~\cite{IKK02}, where 
the asymptotic expansions of the solutions of convection-diffusion equations will be discussed. 

Asymptotic behavior of solutions of nonlinear parabolic equations has been extensively studied 
in many papers by various methods. 
See e.g. \cite{BFW}--\cite{Y02} and references therein. 
Among others, 
Fujigaki and Miyakawa~\cite{FM} studied the large time behavior of the solution~$u$ 
of the Cauchy problem for the incompressible Navier-Stokes equation, 
and gave higher order asymptotic expansions of the solution~$u$ satisfying 
\begin{equation}
\label{eq:1.6}
\sup_{0\le\ell\le N+1}\,\sup_{(x,t)\in{\bf R}^N\times(0,\infty)}\,(1+|x|)^{\ell}(1+t)^{(N+1-\ell)/2}|u(x,t)|
<\infty.
\end{equation}
Their arguments can be also applied to convection-diffusion equations 
(see e.g. \cite{Yamada}--\cite{Y01}). 
On the other hand, in \cite{IK02}
the first and the second authors of this paper  
considered the Cauchy problem for the semilinear heat equation 
\begin{equation}
\label{eq:1.7}
\partial_t u=\Delta u+\lambda|u|^{p-1}u\quad\mbox{in}\quad{\bf R}^N\times(0,\infty),
\end{equation}
where $\lambda\in{\bf R}$ and $p>1+2/N$, 
and gave the precise description of the asymptotic behavior of the solution 
behaving like a multiple of the heat kernel (see also \cite{IIK}). 
Furthermore, in \cite{IK03} 
they extended the results in \cite{IK02}, and established the method of obtaining 
higher order asymptotic expansions of the solutions behaving like a multiple of the heat kernel as $t\to\infty$ 
for general nonlinear heat equations. 
The arguments in \cite{IK03} are applicable to various nonlinear heat equations systematically 
without pointwise decay estimates of the solutions as $|x|\to\infty$, such as \eqref{eq:1.6}. 
\vspace{3pt}

In this paper we improve and generalize 
the arguments in \cite{IIK}, \cite{IK02} and \cite{IK03}, 
and establish the method of obtaining 
the higher order asymptotic expansions of the solutions of nonlinear integral equation~\eqref{eq:1.1} 
behaving like a multiple of the integral kernel $G$ as $t\to\infty$. 
Our arguments are applicable to general nonlinear parabolic equations including \eqref{eq:1.4} and \eqref{eq:1.5}, 
and they can also give some new and sharp decay estimates of the solutions 
even if we focus on the semilinear heat equation~\eqref{eq:1.7} 
(see also Remark~\ref{Remark:1.1}). 
\vspace{5pt}

We introduce some notation. 
For any $k\ge 0$, let $[k]\in{\bf N}\cup\{0\}$ be such that $k-1<[k]\le k$.
For any multi-index $\alpha\in{\bf M}:=({\bf N}\cup\{0\})^N$, 
put 
$$
|\alpha|:=\displaystyle{\sum_{i=1}^N}\alpha_i,\quad
\alpha!:=\prod_{i=1}^N\alpha_i!,\quad
x^\alpha:=\prod_{i=1}^N x_i^{\alpha_i},\quad
\partial_x^\alpha:=
\frac{\partial^{|\alpha|}}{\partial x_1^{\alpha_1}\cdots\partial x_N^{\alpha_N}}.
$$
Let ${\bf M}_k:=\left\{\alpha\in{\bf M}:|\alpha|\le k\right\}$ for $k\ge 0$. 
For any $\alpha=(\alpha_1,\dots,\alpha_N)$, $\beta=(\beta_1,\dots,\beta_N)\in{\bf M}$, 
we say 
$$
\alpha\le\beta
$$
if $\alpha_i\le\beta_i$ for all $i\in\{1,\dots,N\}$.  
For any $1\le r\le\infty$, let $\|\cdot\|_r$ be the usual norm of $L^r:=L^r({\bf R}^N)$. 
For any $k\ge 0$, 
we denote by $|||\,\cdot\,|||_k$ the norm of $L^1_k:=L^1({\bf R}^N,(1+|x|^k)dx)$, that is, 
$$
|||f|||_k:=\int_{{\bf R}^N}|f(x)|(1+|x|^k)dx,\qquad f\in L^1_k. 
$$
For any $\varphi\in L^q({\bf R}^N)$ $(1\le q\le\infty)$, 
we put 
\begin{equation}
\label{eq:1.8}
e^{t{\mathcal L}}\varphi(x):=\int_{{\bf R}^N}G(x-y,t)\varphi(y)dy,\quad x\in{\bf R}^N,\,\,\, t>0. 
\end{equation}
Then, under assumption~$({\bf G})$, we have the following (see also Section~2):
\begin{itemize}
  \item Let $1\le q\le r\le\infty$. Then there exists a constant $C$ such that 
	$$
	\|e^{t{\mathcal L}}\varphi\|_r\le Ct^{-\frac{N}{d}(\frac{1}{q}-\frac{1}{r})}\|\varphi\|_q,\quad t>0,
	$$
	for any $\varphi\in L^q$;
  \item For any $\varphi\in L^q$ with $1\le q\le\infty$, 
	$$
	e^{t{\mathcal L}}\varphi(x)=e^{(t-s){\mathcal L}}[e^{s{\mathcal L}}\varphi](x)
	$$ 
	for all $x\in{\bf R}^N$ and $0<s<t$.
\end{itemize}
Let $0\le k<L$ with $[k]\le\gamma$ and $f\in L^1_k$. 
Put 
$$
g(x,t):=G(x,t+1),\qquad
g_\alpha(x,t):=\frac{(-1)^{|\alpha|}}{\alpha!}\partial_x^\alpha g(x,t)\quad(\alpha\in{\bf M}_\gamma).
$$
Then, for any $t\ge 0$, 
we denote by $P_k(t)f\in L^1_k$ by 
\begin{equation}
\label{eq:1.9}
[P_k(t)f](x):=f(x)-\sum_{|\alpha|\le k}M_\alpha(f,t) g_\alpha(x,t),
\end{equation}
where $M_\alpha(f,t)$ $(|\alpha|\le k)$ are defined inductively (in $\alpha$) by 
\begin{equation}
\label{eq:1.10}
\left\{
\begin{array}{l}
\displaystyle{M_0(f,t):=\int_{{\bf R}^N}f(x)dx}\quad\mbox{if}\quad \alpha=0,\vspace{7pt}\\
\displaystyle{M_\alpha(f,t):=\int_{{\bf R}^N} 
x^\alpha f(x)dx-\sum_{\beta\le\alpha,\beta\not=\alpha}M_\beta(f,t) 
\int_{{\bf R}^N}x^\alpha g_\beta(x,t)dx}\quad\mbox{if}\quad\alpha\not=0. 
\end{array}
\right.
\end{equation}
Then it follows that 
$$
\int_{{\bf R}^N}x^\alpha[P_k(t)f](x)dx=0,\qquad t>0,
$$
for any  $\alpha\in{\bf M}_k$ 
(see Lemma~\ref{Lemma:2.1}~(ii)). 
This is a crucial property of the operator $P_k(t)$ (on $L^1_k$) in our analysis. 
\vspace{5pt}

Now we are ready to state the main results of this paper, 
which give asymptotic expansions of the functions 
\begin{eqnarray*}
e^{t{\cal L}}\varphi(x) \!\!\!& = &\!\!\! \int_{{\bf R}^N}G(x-y,t)\varphi(y)dy,\\
\int_0^te^{(t-s){\cal L}}f(s)ds \!\!\!& = &\!\!\! 
\int_0^t\int_{{\bf R}^N}G(x-y,t-s)f(y,s)dyds,
\end{eqnarray*}
as $t\to\infty$, under suitable integrability conditions on $\varphi$ and $f$. 
\begin{theorem}
\label{Theorem:1.1}
Assume condition~$({\bf G})$ for some $\gamma\in{\bf N}$, $d>0$ and $L>0$. 
Let $0\le K<L$ with $[K]+1\le\gamma$. 
For any $\varphi\in L^1_K$, put 
\begin{equation}
\label{eq:1.11}
v(x,t):=e^{t{\mathcal L}}\varphi(x)-\sum_{|\alpha|\le K}M_\alpha(\varphi,0)g_\alpha(x,t). 
\end{equation}
Then, for any $j\in\{0,\dots,\gamma\}$, $q\in[1,\infty]$ and $\ell\in[0,K]$, 
there exists a constant $C$ independent of $\varphi\in L^1_K$ such that 
\begin{equation}
t^{\frac{N}{d}(1-\frac{1}{q})+\frac{j}{d}}\|\nabla^j v(t)\|_q+t^{\frac{j}{d}}(1+t)^{-\frac{\ell}{d}}|||\nabla^j v(t)|||_\ell
\le C(1+t)^{-\frac{K}{d}}|||\varphi|||_K
\label{eq:1.12}
\end{equation}
for all $t>0$. 
Furthermore, for any $\varphi\in L^1_K$, 
\begin{equation}
\lim_{t\to\infty}t^{\frac{K}{d}}
\biggr[t^{\frac{N}{d}(1-\frac{1}{q})+\frac{j}{d}}\|\nabla^j v(t)\|_q+t^{-\frac{\ell}{d}+\frac{j}{d}}|||\nabla^j v(t)|||_\ell\biggr]=0.
\label{eq:1.13}
\end{equation}
\end{theorem}
\begin{theorem}
\label{Theorem:1.2}
Assume condition~$({\bf G})$ for some $\gamma\in{\bf N}$, $d>0$ and $L>0$. 
Let $0\le K<L$ with $[K]+1\le\gamma$ and $1\le q\le\infty$. 
Let $f$ be a measurable function in ${\bf R}^N\times(0,\infty)$ such that 
\begin{equation}
\label{eq:1.14}
E_{K,q}[f](t):=(1+t)^{\frac{K}{d}}\left[t^{\frac{N}{d}(1-\frac{1}{q})}\|f(t)\|_q+\|f(t)\|_1\right]+|||f(t)|||_K\in L^\infty(0,T)
\end{equation}
for any $T>0$. 
Then the following holds: 
\vspace{3pt}
\newline
{\rm (i)} For any $\alpha\in{\bf M}_K$, there exists a constant $C_1$ such that
\begin{equation}
\label{eq:1.15}
|M_\alpha(f(t),t)|\le C_1(1+t)^{-\frac{K-|\alpha|}{d}}E_{K,q}[f](t)
\end{equation}
for almost all $t>0$;
\vspace{3pt}
\newline
{\rm (ii)} 
Put 
\begin{eqnarray*}
 R_K[f](t) \!\!\!& := &\!\!\!
 \int_0^t e^{(t-s){\cal L}}P_K(s)f(s)ds\\
 \!\!\! & \,\,= &\!\!\!
\int_0^t e^{(t-s){\mathcal L}}f(s)ds-\sum_{|\alpha|\le K}\left[\int_0^t M_\alpha(f(s),s)ds\right]g_\alpha(t).
\end{eqnarray*}
Let $j\in\{0,\dots,\gamma\}$ with $j<d$ and $T_0>0$. 
Then there exists a constant $C_2$ such that, 
for any $\epsilon>0$ and $T\ge T_0$, 
\begin{eqnarray}
 & & t^{\frac{N}{d}(1-\frac{1}{q})}\|\nabla^jR_K[f](t)\|_q+t^{-\frac{\ell}{d}}|||\nabla^j R_K[f](t)|||_\ell\notag\\
 & & \qquad\quad
 \le \epsilon t^{-\frac{K+j}{d}}
+C_2t^{-\frac{K}{d}}\int_T^t (t-s)^{-\frac{j}{d}}E_{K,q}[f](s)ds
\label{eq:1.16}
\end{eqnarray}
for all sufficiently large $t>0$. 
In particular, if 
$$
\int_0^\infty E_{K,q}[f](s)ds<\infty,
$$
then 
\begin{equation}
\label{eq:1.17}
\lim_{t\to\infty}
t^{\frac{K}{d}}
\biggr[t^{\frac{N}{d}(1-\frac{1}{q})}\|R_K[f](t)\|_q+t^{-\frac{\ell}{d}}|||R_K[f](t)|||_\ell\biggr]=0.
\end{equation}
\end{theorem}
By Theorems~\ref{Theorem:1.1} and \ref{Theorem:1.2} we can give  
decay estimates of the distance in $L^q$ and $L^1_\ell$ $(0\le \ell\le K)$ 
from the solution of \eqref{eq:1.1} to its asymptotic expansion  
\begin{equation}
\label{eq:1.18}
\sum_{|\alpha|\le K}\left[M_\alpha(\varphi,0)+\int_0^t M_\alpha(f(s),s)ds\right]g_\alpha(t). 
\end{equation}
The higher order asymptotic expansions of the solutions depend on the nonlinearity of $f$ 
and are discussed in Sections~4 and 5. 
\begin{remark}
\label{Remark:1.1}
Let $G$ be the heat kernel, that is, 
$$
G(x,t):=(4\pi t)^{-\frac{N}{2}}\exp\left(-\frac{|x|^2}{4t}\right).
$$
Let $\varphi\in L^1_K$ for some $K\ge 0$, 
and define a function $v$ by \eqref{eq:1.11}. 
In {\rm\cite{IIK}} the authors proved that, 
for any $1\le q\le\infty$ and $0\le\ell\le K$, 
\begin{equation}
\label{eq:1.19}
t^{\frac{N}{2}(1-\frac{1}{q})}\|v(t)\|_q
=\left\{
\begin{array}{ll}
O(t^{-\frac{K}{2}}) & \mbox{if}\quad K>[K],\\
o(t^{-\frac{K}{2}}) & \mbox{if}\quad K=[K],
\end{array}
\right.
\quad
t^{-\frac{\ell}{2}}|||v(t)|||_\ell=O(t^{-\frac{K}{2}+\sigma}),
\end{equation}
as $t\to\infty$, for any $\sigma>0$. 
This is one of the main ingredients of the asymptotic analysis in {\rm\cite{IIK}}, {\rm\cite{IK02}} and {\rm\cite{IK03}} 
for parabolic equations. 

On the other hand, since the heat kernel satisfies condition~$({\bf G})$ 
for any $\gamma\in{\bf N}$ and $L>0$ with $d=2$, 
Theorem~{\rm\ref{Theorem:1.1}} gives better decay estimates of $v$ than \eqref{eq:1.19}, 
and enables us to improve the asymptotic analysis in {\rm\cite{IIK}}, {\rm\cite{IK02}} and {\rm\cite{IK03}}. 
See Sections~{\rm 5} and {\rm 6}.
\end{remark}

The rest of this paper is organized as follows. 
Section~2 presents some preliminaries on $e^{t{\mathcal L}}\varphi$ and $M_\alpha(f,t)$. 
In Section~3 we improve the argument in \cite{IIK}, 
and study the asymptotic expansion of $e^{t{\mathcal L}}\varphi$. 
This enables us to prove Theorem~\ref{Theorem:1.1}. 
Section~4 is devoted to the proof of Theorem~\ref{Theorem:1.2} 
by using the arguments in the previous sections. 
In Section~5 we study the asymptotic behavior of solutions 
of integral equations with power nonlinearity. 
In Section~6 we apply our arguments to semilinear parabolic equations~\eqref{eq:1.4} and \eqref{eq:1.5},  
and show the validity of our arguments.
\section{Preliminaries}
In this section 
we prove some preliminary results on $e^{t{\mathcal L}}\varphi$ and $M_\alpha(f,t)$. 
In what follows, 
for any two nonnegative functions 
$f_1$ and $f_2$ in a subset $D$ of $[0,\infty)$, 
we say 
$$
f_1(t)\preceq f_2(t),\qquad t\in D
$$ 
if there exists a positive constant $C$ such that 
$f_1(t)\le Cf_2(t)$ for all $t\in D$. 
In addition, we say  
$$
f_1(t)\asymp f_2(t),\qquad t\in D
$$ 
if $f_1(t)\preceq  f_2(t)$ and $f_2(t)\preceq f_1(t)$ for all $t\in D$. 
\vspace{5pt}

We first state some properties on the kernel $G$, which immediately follow 
from condition~$({\bf G})$ (see also \cite{IKK01}):
\begin{itemize}
  \item[{\rm (i)}] 
	$\displaystyle{\int_{{\bf R}^N}G(x,t)dx=1}$ for any $t>0$;
  \item[{\rm (ii)}] 
  	For any $\alpha\in{\bf M}_\gamma$, 
	\begin{equation}
	\label{eq:2.1}
	|\partial_x^\alpha G(x,t)|\preceq t^{-\frac{N}{d}-\frac{|\alpha|}{d}}
	\left(1+\frac{|x|}{t^{1/d}}\right)^{-(N+L+|\alpha|)},
	\quad(x,t)\in{\bf R}^N\times(0,\infty);
		\end{equation}
  \item[{\rm (iii)}] 
	For any $1\le r\le\infty$, $\alpha\in{\bf M}_\gamma$ and $\ell\in[0,L+|\alpha|)$, 
	\begin{equation}
	\sup_{t>0}\,\left[t^{\frac{N}{d}(1-\frac{1}{r})+\frac{|\alpha|}{d}}\|\partial_x^\alpha G(t)\|_r
	+(1+t)^{-\frac{\ell}{d}}t^{\frac{|\alpha|}{d}}|||\partial_x^\alpha G(t)|||_\ell\right]<\infty.
	\label{eq:2.2}	
	\end{equation}
\end{itemize}
Furthermore, applying the Young inequality to \eqref{eq:1.8} with the aid of property~(iii), 
for any $1\le r\le q\le\infty$ and $j\in\{0,\dots,\gamma\}$, 
we can find a constant $C$ such that 
\begin{equation}
\label{eq:2.3}
\|\nabla^j e^{t{\mathcal L}}\varphi\|_q\preceq t^{-\frac{N}{d}(\frac{1}{r}-\frac{1}{q})-\frac{j}{d}}\|\varphi\|_r,
\quad t>0, 
\end{equation}
for all $\varphi\in L^r$ . 
In particular, for $r=q$, by property~(i) we have 
\begin{equation}
\label{eq:2.4}
\|e^{t{\mathcal L}}\varphi\|_q\le\|\varphi\|_q,
\quad t>0. 
\end{equation}
Next we state a lemma on $M_\alpha(f,t)$ and the operator $P_k(t)$. 
\begin{lemma}
\label{Lemma:2.1}
Assume condition~$({\bf G})$ for some $\gamma\in{\bf N}$, $d>0$ and $L>0$.  
For any $f$, $g\in L^1_k$ with $0\le k<L$ and $[k]\le\gamma$, 
the following holds: 
\vspace{3pt}
\newline
{\rm (i)} 
For any $a,b\in{\bf R}$ and $\alpha\in{\bf M}_k$, 
$$
M_\alpha(af+bg,t)=aM_\alpha(f,t)+bM_\alpha(g,t),
\qquad t\ge 0;
$$
{\rm (ii)} 
For any $\alpha\in {\bf M}_k$, 
$$
\int_{{\bf R}^N}x^\alpha[P_k(t)f](x)dx=0,\qquad t\ge 0;
$$
{\rm (iii)} 
Assume that there exists constants $\{c_\alpha\}_{\alpha\in{\bf M}_k}$ such that 
$$
\int_{{\bf R}^N}x^\beta\biggr(f-\sum_{|\alpha|\le k}c_\alpha g_\alpha(x,t)\biggr)dx=0,\qquad
\beta\in{\bf M}_k,
$$
for some $t\ge 0$. Then 
$$
c_\alpha=M_\alpha(f,t),\qquad \alpha\in{\bf M}_k;
$$
{\rm (iv)} For any $t\ge 0$,  
$$
M_\alpha(e^{t{\mathcal L}}f,t)=M_\alpha(f,0),\qquad \alpha\in{\bf M}_k;
$$
{\rm (v)} Let $f\in L^1_k$ be such that 
\begin{equation}
\label{eq:2.5}
\int_{{\bf R}^N}x^\beta f(x)dx=0,\qquad \beta\in {\bf M}_k. 
\end{equation}
Then 
$$
\int_{{\bf R}^N}x^\beta e^{t{\mathcal L}}f(x)dx=0,\qquad \beta\in {\bf M}_k,
$$
for all $t\ge 0$. 
\end{lemma}
{\bf Proof.}
Assertion~(i) immediately follows from \eqref{eq:1.10}. 
We prove assertion~(ii). 
For any $f\in L^1_k$ and $\alpha\in{\bf M}_k$, 
since 
$$
\int_{{\bf R}^N}x^\alpha g_\beta(x,t)dx=0
\quad\mbox{if not $\beta\le\alpha$}, 
$$
by \eqref{eq:1.9} and \eqref{eq:1.10} we have 
\begin{eqnarray*}
 & & \int_{{\bf R}^N}x^\alpha[P_k(t)f](x)dx 
=\int_{{\bf R}^N}x^\alpha f(x)dx-\sum_{\beta\le\alpha}M_\beta(f,t)\int_{{\bf R}^N}x^\alpha g_\beta(x,t)dx\\
 & & \qquad\qquad
=\int_{{\bf R}^N}x^\alpha f(x)dx-M_\alpha(f,t)
-\sum_{\beta\le\alpha,\alpha\not=\beta}M_\beta(f,t)\int_{{\bf R}^N}x^\alpha g_\beta(x,t)dx=0
\end{eqnarray*}
for $t\ge 0$. This implies assertion~(ii). 
Similarly assertion~(iii) follows inductively in $\alpha$. 

We prove assertion~(iv). 
For any $f\in L^1_k$, put 
\begin{equation}
\label{eq:2.6}
w(x,t):=e^{t{\mathcal L}}f(x)-\sum_{|\alpha|\le k}M_\alpha(f,0)g_\alpha(x,t). 
\end{equation}
Since 
\begin{eqnarray*}
 & & [e^{t{\mathcal L}}g_\alpha(0)](x)
=\frac{(-1)^{|\alpha|}}{\alpha!}\int_{{\bf R}^N}G(x-y,t)\partial_y^\alpha G(y,1)dy\\
 & & \qquad
=\frac{(-1)^{|\alpha|}}{\alpha!}\partial_x^\alpha\int_{{\bf R}^N}G(x-y,t)G(y,1)dy
=\frac{(-1)^{|\alpha|}}{\alpha!}\partial_x^\alpha G(x,t+1)
=g_\alpha(x,t),
\end{eqnarray*}
we have 
$$
w(x,t)=[e^{t{\mathcal L}}w(0)](x).
$$
On the other hand, it follows from assertion~(ii) that 
$$
\int_{{\bf R}^N} x^\beta w(x,0)dx
=\int_{{\bf R}^N}x^\beta[P_k(0)f](x)dx=0,\qquad \beta\in{\bf M}_k. 
$$
Therefore, by the Fubini theorem and the binomial theorem we have 
\begin{eqnarray*}
\int_{{\bf R}^N} x^\beta w(x,t)dx
\!\!\! & = &\!\!\!
\int_{{\bf R}^N}x^\beta\left(\int_{{\bf R}^N}G(x-y,t)w(y,0)dy\right)dx\\
 \!\!\! & = &\!\!\!
 \int_{{\bf R}^N}\left(\int_{{\bf R}^N}(x+y)^\beta G(x,t)dx\right)w(y,0)dy\\
  \!\!\! & = &\!\!\!
  \sum_{\alpha\le\beta}C_\alpha(t)\int_{{\bf R}^N}y^\alpha w(y,0)dy=0
\end{eqnarray*}
for $\beta\in{\bf M}_k$, 
where $\{C_\alpha(t)\}$ are constants depending on $t$. 
Then assertion~(iv) follows from assertion~(iii) and \eqref{eq:2.6}. 

It remains to prove assertion~(v). 
Let $f\in L^1_k$, and assume \eqref{eq:2.5}. 
By \eqref{eq:1.10} we obtain inductively 
$$
M_\alpha(f,0)=0,\qquad \alpha\in{\bf M}_k. 
$$
This together with assertion~(iv) implies that 
\begin{eqnarray*}
[P_k(t)e^{t{\mathcal L}}f](x)
\!\!\! & = &\!\!\!
e^{t{\mathcal L}}f(x)-\sum_{|\alpha|\le k}M_\alpha(e^{t{\mathcal L}}f,t)g_\alpha(x,t)\\
\!\!\! & = &\!\!\!
e^{t{\mathcal L}}f(x)-\sum_{|\alpha|\le k}M_\alpha(f,0)g_\alpha(x,t)
=e^{t{\mathcal L}}f(x).
\end{eqnarray*}
Then assertion~(v) follows from assertion~(ii), 
and the proof of Lemma~\ref{Lemma:2.1} is complete. 
$\Box$\vspace{5pt}

At the end of this section, 
we prove a lemma on the functions $P_K(t)f(t)$ and $E_K[f](t)$.
\begin{lemma}
\label{Lemma:2.2}
Assume the same conditions as in Theorem~{\rm\ref{Theorem:1.2}}. 
Let $1\le r\le q\le\infty$ and $0\le\ell\le K$. 
Then 
\begin{eqnarray}
 & & \|f(t)\|_r\le t^{-\frac{N}{d}(1-\frac{1}{r})}(1+t)^{-\frac{K}{d}}E_{K,q}[f](t),
  \label{eq:2.7}\\
 & & |||f(t)|||_\ell\preceq (1+t)^{-\frac{K}{d}+\frac{\ell}{d}}E_{K,q}[f](t),
 \label{eq:2.8}\\
 & & |M_\alpha(f(t),t)|\preceq 
 (1+t)^{-\frac{K}{d}+\frac{|\alpha|}{d}}E_{K,q}[f](t),\quad \alpha\in{\bf M}_K,
  \label{eq:2.9}
\end{eqnarray}
for almost all $t>0$. 
Furthermore, 
\begin{equation}
t^{\frac{N}{d}(1-\frac{1}{q})}\|P_K(t)f(t)\|_q+(1+t)^{-\frac{\ell}{d}}|||P_K(t)f(t)|||_\ell
\preceq (1+t)^{-\frac{K}{d}}E_{K,q}[f](t)
\label{eq:2.10}
\end{equation}
for almost all $t>0$. 
\end{lemma}
{\bf Proof.} 
For $1\le r\le q$, by the H\"older inequality and \eqref{eq:1.14} we have 
$$
\|f(t)\|_r\le\|f(t)\|_1^{1-\frac{\eta}{r}}\|f(t)\|_q^{\frac{\eta}{r}}
\le t^{-\frac{N}{d}(1-\frac{1}{r})}(1+t)^{-\frac{K}{d}}E_{K,q}[f](t)
$$
for almost all $t>0$, where $\eta:=(r-1)/(1-1/q)$, and we obtain \eqref{eq:2.7}. 
For any $0\le \ell\le K$, since 
$$
\left(\frac{1+|x|}{(1+t)^{1/d}}\right)^\ell \preceq 1+\left(\frac{1+|x|}{(1+t)^{1/d}}\right)^{K},
$$
we get 
\begin{eqnarray*}
 & & |||f(t)|||_\ell\preceq
(1+t)^{\frac{\ell}{d}}\int_{{\bf R}^N}\left(\frac{1+|x|}{(1+t)^{1/d}}\right)^\ell |f(x,t)|dx\\
 & & \qquad\quad
 \preceq(1+t)^{\frac{\ell}{d}}\int_{{\bf R}^N}\biggr[1+\left(\frac{1+|x|}{(1+t)^{1/d}}\right)^{K}\biggr]|f(x,t)|dx\\
 & & \qquad\quad
 \preceq(1+t)^{\frac{\ell}{d}}\biggr[\|f(t)\|_1+(1+t)^{-\frac{K}{d}}|||f(t)|||_K\biggr]
\le (1+t)^{-\frac{K}{d}+\frac{\ell}{d}}E_{K,q}[f](t)
\end{eqnarray*}
for almost all $t>0$. This implies \eqref{eq:2.8}. 

The proof of \eqref{eq:2.9} is by induction in $\alpha\in M_K$. 
For $\alpha=0$, by \eqref{eq:1.10} and \eqref{eq:2.7} we have 
$$
|M_\alpha(f(t),t)|=\left|\int_{{\bf R}^N}f(x,t)dx\right|\le \|f(t)\|_1\le (1+t)^{-\frac{K}{d}}E_{K,q}[f](t)
$$
for almost all $t>0$, and inequality~\eqref{eq:2.9} holds for $\alpha=0$.  
Assume that inequality~\eqref{eq:2.9} holds for $\alpha\in{\bf M}_n$ for some $n\in\{0,\dots,[K]-1\}$. 
Then, for any $\alpha\in{\bf M}_{n+1}\setminus{\bf M}_n$, 
it follows from \eqref{eq:1.10}, \eqref{eq:2.2} and \eqref{eq:2.8} that 
\begin{eqnarray*}
|M_\alpha(f,t)| & \le & \left|\int_{{\bf R}^N}x^\alpha f(x, t)dx\right|
+\sum_{\beta\le\alpha,\,\beta\not=\alpha}|M_\beta(f(t),t)|\left|\int_{{\bf R}^N}x^\alpha g_\beta(x,t)dx\right|\\
 & \preceq & |||f(t)|||_{|\alpha|}
 +\sum_{\beta\le\alpha,\beta\not=\alpha}(1+t)^{-\frac{K}{d}+\frac{|\beta|}{d}}E_{K,q}[f](t)\cdot(1+t)^{\frac{|\alpha|}{d}-\frac{|\beta|}{d}}\\
 & \preceq & (1+t)^{-\frac{K}{d}+\frac{|\alpha|}{d}}E_{K,q}(t)
\end{eqnarray*}
for almost all $t>0$, and \eqref{eq:2.9} holds. 
Thus inequality \eqref{eq:2.9} holds for all $\alpha\in{\bf M}_K$. 
Furthermore, 
by \eqref{eq:1.9}, \eqref{eq:2.2} and \eqref{eq:2.9} 
we obtain 
\begin{eqnarray*}
 & & t^{\frac{N}{d}(1-\frac{1}{q})}\|P_K(t)f(t)\|_q+(1+t)^{-\frac{\ell}{d}}|||P_K(t)f(t)|||_\ell\\
 & & \le t^{\frac{N}{d}(1-\frac{1}{q})}\|f(t)\|_q+(1+t)^{-\frac{\ell}{d}}|||f(t)|||_\ell\\
 & & \qquad\quad
 +\sum_{|\alpha|\le K}|M_\alpha(f(t),t)|
 \left[t^{\frac{N}{d}(1-\frac{1}{q})}\|g_\alpha(t)\|_q+(1+t)^{-\frac{\ell}{d}}|||g_\alpha(t)|||_\ell\right]\\
 & & \preceq (1+t)^{-\frac{K}{d}}E_{K,q}[f](t)
\end{eqnarray*}
for all $t>0$. This implies \eqref{eq:2.10}, 
and the proof of Lemma~\ref{Lemma:2.2} is complete. 
$\Box$
\section{Proof of Theorem~\ref{Theorem:1.1}}
In this section we prove the following proposition on 
the decay estimates of $e^{t{\mathcal L}}\varphi$. 
Proposition~\ref{Proposition:3.1} is one of the main ingredients of this paper 
and improves \cite[Lemmas~2.2 and 2.5]{IIK}. 
Theorem~\ref{Theorem:1.1} follows from Lemma~\ref{Lemma:2.1} and Proposition~\ref{Proposition:3.1}. 
\begin{proposition}
\label{Proposition:3.1}
Assume condition~$({\bf G})$ for some $\gamma\in{\bf N}$, $d>0$ and $L>0$. 
Let $0\le k<L$ with $[k]+1\le\gamma$ and $j\in\{0,\dots,\gamma\}$. 
\vspace{3pt}
\newline
{\rm (i)} 
For any $\ell\in[0,k]$, 
there exists a constant $C_1$ such that 
\begin{equation}
 \int_{{\bf R}^N}|x|^\ell|\nabla^j e^{t{\mathcal L}}\varphi(x)|dx
 \le C_1t^{-\frac{j-\ell}{d}}\int_{{\bf R}^N}|\varphi(x)|dx+C_1t^{-\frac{j}{d}}\int_{{\bf R}^N}|x|^\ell|\varphi(x)|dx
\label{eq:3.1}
\end{equation}
for all $t>0$ and $\varphi\in L^1_k$. 
\vspace{3pt}
\newline
{\rm (ii)}
For any $\ell\in[0,k]$, 
there exists a constant $C_2$ such that 
\begin{equation}
\label{eq:3.2}
\int_{{\bf R}^N}|x|^{\ell}|e^{t{\cal L}}\varphi(x)|dx\le C_2t^{-\frac{k-\ell}{d}}\int_{{\bf R}^N}|x|^k|\varphi(x)|dx,
\quad t>0,
\end{equation}
for all $\varphi\in L^1_k$ satisfying 
\begin{equation}
\label{eq:3.3}
\int_{{\bf R}^N}x^\alpha\varphi(x)dx=0,
\quad \alpha\in{\bf M}_k.
\end{equation}
{\rm (iii)} 
For any $\ell\in[0,k]$,  
\begin{equation}
\label{eq:3.4}
\lim_{t\to\infty}t^{\frac{k-\ell}{d}}\int_{{\bf R}^N}|x|^k|e^{t{\cal L}}\varphi(x)|dx=0
\end{equation}
for $\varphi\in L^1_k$ satisfying \eqref{eq:3.3}.
\end{proposition}
{\bf Proof.} 
Let $\varphi\in L^1_k$ and $0\le \ell\le k$. 
By \eqref{eq:2.2} 
we have 
\begin{eqnarray}
 & & \int_{{\bf R}^N}|x|^\ell|\nabla^jG(x-y,t)|dx
=\int_{{\bf R}^N}|x+y|^\ell|\nabla^jG(x,t)|dx\notag\\
 & & \qquad\qquad
 \preceq\int_{{\bf R}^N}(|x|^\ell+|y|^\ell)|\nabla^jG(x,t)|dx
\preceq t^{-\frac{j-\ell}{d}}+t^{-\frac{j}{d}}|y|^\ell
\label{eq:3.5}
\end{eqnarray}
for all $y\in{\bf R}^N$ and $t>0$. 
This implies  
\begin{eqnarray*}
 \int_{{\bf R}^N}|x|^\ell |\nabla^j e^{t{\mathcal L}}\varphi(x)|dx
 \!\!\! & \le &\!\!\!\int_{{\bf R}^N}
 \left(\int_{{\bf R}^N}|x|^\ell |\nabla^jG(x-y,t)|dx\right)|\varphi(y)|dy\\
 \!\!\!& \preceq &\!\!\! t^{-\frac{j-\ell}{d}}\|\varphi\|_1+t^{-\frac{j}{d}}\int_{{\bf R}^N}|x|^\ell |\varphi(x)|dx
\end{eqnarray*}
for all $t>0$, and assertion~(i) follows.  

In order to prove assertions~(ii) and (iii), we assume \eqref{eq:3.3}. 
Let $R>1$, and put 
\begin{eqnarray*}
 & & \varphi_1(x):=\varphi(x)\chi_{\{|x|<R^{1/d}\}}(x),\qquad
 \varphi_2(x):=\varphi(x)\chi_{\{|x|\ge R^{1/d}\}}(x),
 \vspace{3pt}\\
 & & \psi_R(x):=[P_k(R)\varphi_1](x)
 =\varphi_1(x)-\tilde{\psi}_R(x),
\,\,\,\,
\tilde{\psi}_R(x):=\sum_{|\beta|\le k}M_\beta(\varphi_1,R)g_\beta(x,R).
\end{eqnarray*}
We prove that, for any $\beta\in{\bf M}_k$,    
there exists a constant $C_1$ such that 
\begin{equation}
\label{eq:3.6}
|M_\beta(\varphi_1,R)|\le C_1R^{-\frac{k-|\beta|}{d}}\int_{\{|x|\ge R^{1/d}\}}|x|^k|\varphi(x)|dx,
\quad R>1. 
\end{equation}
For any $\beta\in{\bf M}_k$, by \eqref{eq:3.3} we obtain 
\begin{eqnarray}
 & & \left|\int_{{\bf R}^N}x^\beta\varphi_1(x)dx\right|
=\left|\int_{\{|x|<R^{1/d}\}}x^\beta\varphi(x)dx\right|
=\left|\int_{\{|x|\ge R^{1/d}\}}x^\beta\varphi(x)dx\right|\notag\\
 & & \qquad\quad
 \le\int_{\{|x|\ge R^{1/d}\}}|x|^{|\beta|}|\varphi(x)|dx
 \le R^{-\frac{k-|\beta|}{d}}\int_{\{|x|\ge R^{1/d}\}}|x|^k|\varphi(x)|dx.
 \label{eq:3.7}
\end{eqnarray}
This implies \eqref{eq:3.6} for $\beta\in{\bf M}_0$, 
and \eqref{eq:3.6} holds in the case $0\le k<1$. 
In the case $k\ge 1$, 
we assume that inequality~\eqref{eq:3.6} holds for all $\beta\in{\bf M}_n$, 
where $n\in\{0,\dots,[k]-1\}$. 
Then, for any $\beta\in{\bf M}_k$ with $|\beta|=n+1$, 
by \eqref{eq:1.10}, \eqref{eq:2.2} and \eqref{eq:3.7} 
we have 
\begin{eqnarray*}
|M_\beta(\varphi_1,R)| 
\!\!\!& \le &\!\!\! 
\left|\int_{{\bf R}^N}x^\beta\varphi_1(x)dx\right|
+\sum_{\alpha\le\beta,\,\alpha\not=\beta}
|M_\alpha(\varphi_1,R)|\left|\int_{{\bf R}^N}x^\beta g_\alpha(x,R)dx\right|\\
 \!\!\! & \le &\!\!\!  
 \left|\int_{{\bf R}^N}x^\beta\varphi_1(x)dx\right|
 +C_2\sum_{\alpha\le\beta,\,\alpha\not=\beta}
 (1+R)^{\frac{|\beta|-|\alpha|}{d}}|M_\alpha(\varphi_1,R)|\\
 \!\!\!& \le &\!\!\! 
 C_3R^{-\frac{k-|\beta|}{d}}\int_{\{|x|\ge R^{1/d}\}}|x|^k|\varphi(x)|dx,
 \qquad R>1,
\end{eqnarray*}
for some constants $C_2$ and $C_3$. 
This implies \eqref{eq:3.6} for any $\beta\in{\bf M}_k$ with $|\beta|={n+1}$. 
Thus it follows \eqref{eq:3.6} by induction. 

Let $0\le\ell\le k$. 
We apply a similar argument to \eqref{eq:3.7} with the aid of assertion~(i), 
and obtain 
\begin{eqnarray}
 & & \int_{{\bf R}^N}|x|^\ell |e^{t{\mathcal L}}\varphi_2(x)|dx 
\preceq t^{\frac{\ell}{d}}\int_{{\bf R}^N}|\varphi_2(x)|dx
+\int_{{\bf R}^N}|x|^\ell |\varphi_2(x)|dx\notag\\
 & & \qquad\qquad
\preceq (t^{\frac{\ell}{d}}R^{-\frac{k}{d}}+R^{-\frac{k-\ell}{d}})\int_{\{|x|\ge R^{1/d}\}}|x|^k|\varphi(x)|dx\notag\\
 & & \qquad\qquad
\preceq t^{-\frac{k-\ell}{d}}\left[(t^{-1}R)^{-\frac{k}{d}}+(t^{-1}R)^{-\frac{k-\ell}{d}}\right]
\int_{\{|x|\ge R^{1/d}\}}|x|^k|\varphi(x)|dx
\label{eq:3.8}
\end{eqnarray}
for all $t>0$ and $R>1$. 
On the other hand,  
by \eqref{eq:2.2} and \eqref{eq:3.6} 
we have 
\begin{eqnarray}
 & & \int_{{\bf R}^N}|x|^\ell |e^{t{\mathcal L}}\tilde{\psi}_R(x)|dx
\le\sum_{|\beta|\le k}|M_\beta(\varphi_1,R)|\int_{{\bf R}^N}|x|^\ell|g_\beta(x,t+R)|dx\nonumber\\
 & & \qquad\qquad
 \preceq R^{-\frac{k}{d}}\sum_{|\beta|\le k}
\frac{R^{\frac{|\beta|}{d}}}{(1+t+R)^{\frac{|\beta|-\ell}{d}}}
 \int_{\{|x|\ge R^{1/d}\}}|x|^k|\varphi(x)|dx\nonumber\\
 & & \qquad\qquad
 \preceq R^{-\frac{k}{d}}(t+R)^{\frac{\ell}{d}}
 \int_{\{|x|\ge R^{1/d}\}}|x|^k|\varphi(x)|dx\label{eq:3.9}
\end{eqnarray}
for all $t\ge 0$ and $R>1$. 
In particular, by \eqref{eq:3.9} we have 
\begin{eqnarray}
\label{eq:3.10}
 & & \int_{{\bf R}^N}|x|^\ell |\tilde{\psi}_R(x)|dx
\preceq R^{-\frac{k-\ell}{d}}\int_{\{|x|\ge R^{1/d}\}}|x|^k|\varphi(x)|dx,\\
 & & \label{eq:3.11}
\int_{{\bf R}^N}|x|^k |\psi_R(x)|dx\le 
\int_{{\bf R}^N}|x|^k[|\varphi_1(x)|+|\tilde{\psi}_R(x)|]dx
\preceq\int_{{\bf R}^N}|x|^k|\varphi(x)|dx,
\end{eqnarray}
for all $R>1$. 

Put
$$
G_k(x,y,t):=
G(x-y,t)-\sum_{|\alpha|\le k}\frac{(-1)^{|\alpha|}}{\alpha!}\partial_x^\alpha G(x,t)y^\alpha. 
$$
It follows from Lemma~\ref{Lemma:2.1}~(ii) that 
$$
\int_{{\bf R}^N}y^\alpha\psi_R(y)dy=0,\qquad \alpha\in{\bf M}_k. 
$$
This implies that 
\begin{eqnarray}
e^{t{\mathcal L}}\psi_R(x)
\!\!\! & = & \!\!\!\! \int_{{\bf R}^N}G_k(x,y,t)\psi_R(y)dy\nonumber\\
\!\!\! & = &\!\!\! \int_{\{|y|<R^{1/d}\}}G_k(x,y,t)\psi_R(y)dy
-\int_{\{|y|\ge R^{1/d}\}}G_k(x,y,t)\tilde{\psi}_R(y)dy\nonumber\\
\!\!\! & =: &\!\!\! I_1(x,t)-I_2(x,t).\label{eq:3.12}
\end{eqnarray}
By the mean value theorem, for any $y\in{\bf R}^N$, 
we can find $\tilde{y}\in{\bf R}^N$ with $|\tilde{y}|\le|y|$ such that 
$$
|G_k(x,y,t)|\le C_4|\nabla_x^{[k]+1}G(x-\tilde{y},t)||y|^{[k]+1},
$$
where $C_4$ is a constant independent of $y$. 
Then,  by \eqref{eq:2.2}, \eqref{eq:3.11} and \eqref{eq:3.12} we have 
\begin{eqnarray}
 & & \int_{{\bf R}^N}|x|^\ell|I_1(x,t)|dx\notag\\
 & & \preceq\int_{\{|y|<R^{1/d}\}}
 \biggr(\int_{{\bf R}^N}|x|^\ell |\nabla_x^{[k]+1}G(x-\tilde{y},t)|dx\biggr)|y|^{[k]+1}|\psi_R(y)|dy\notag\\
 & & \preceq\int_{\{|y|<R^{1/d}\}}
 \biggr(\int_{{\bf R}^N}(|x|^\ell+|y|^\ell)|\nabla_x^{[k]+1}G(x,t)|dx\biggr)|y|^{[k]+1}|\psi_R(y)|dy\notag\\
 & & \preceq\int_{\{|y|<R^{1/d}\}}
 \biggr[t^{-\frac{[k]+1-\ell}{d}}+|y|^\ell t^{-\frac{[k]+1}{d}}\biggr]|y|^{[k]+1}|\psi_R(y)|dy\notag\\
 & & \preceq\biggr[t^{-\frac{[k]+1-\ell}{d}}R^{\frac{[k]+1-k}{d}}
 +t^{-\frac{[k]+1}{d}}R^{\frac{[k]+1+\ell-k}{d}}\biggr]\int_{{\bf R}^N}|y|^k|\psi_R(y)|dy\notag\\
 & & \preceq t^{-\frac{k-\ell}{d}}
 \left[(t^{-1}R)^{\frac{[k]+1-k}{d}}+(t^{-1}R)^{\frac{[k]+1+\ell-k}{d}}\right]\int_{{\bf R}^N}|y|^k|\varphi(y)|dy
\label{eq:3.13}
\end{eqnarray}
for all $t>0$ and $R>1$. On the other hand, since 
$$
I_2(x,t)=\int_{\{|y|\ge R^{1/d}\}}G(x-y,t)\tilde{\psi}_R(y)dy
-\sum_{|\alpha|\le k}\frac{(-1)^{|\alpha|}}{\alpha!}\partial_x^\alpha G(x,t)
\int_{\{|y|\ge R^{1/d}\}}y^\alpha\tilde{\psi}_R(y)dy,
$$
by \eqref{eq:2.2}, \eqref{eq:3.5} and \eqref{eq:3.10} 
we have 
\begin{eqnarray}
 & & \int_{{\bf R}^N}|x|^\ell|I_2(x,t)|dx\notag\\
 & & \le
 \int_{\{|y|\ge R^{1/d}\}}\left(\int_{{\bf R}^N}|x|^\ell |G(x-y,t)|dx\right)|\tilde{\psi}_R(y)|dy\notag\\
 & & \qquad\qquad
 +\sum_{|\alpha|\le k}\frac{1}{\alpha!}\int_{{\bf R}^N}|x|^\ell |\partial_x^\alpha G(x,t)|dx
 \int_{\{|y|\ge R^{1/d}\}}|y|^{|\alpha|}|\tilde{\psi}_R(y)|dy\notag\\
 & & \preceq
 \int_{\{|y|\ge R^{1/d}\}}(t^{\frac{\ell}{d}}+|y|^\ell)|\tilde{\psi}_R(y)|dy\notag\\
  & & \qquad\qquad\qquad
 +\sum_{|\alpha|\le k}t^{-\frac{|\alpha|-\ell}{d}}R^{-\frac{k-|\alpha|}{d}}\int_{\{|y|\ge R^{1/d}\}}|y|^k|\varphi(y)|dy\notag\\
 & & \preceq
 t^{-\frac{k-\ell}{d}}\biggr[(t^{-1}R)^{-\frac{k}{d}}+(t^{-1}R)^{-\frac{k-\ell}{d}}+\sum_{|\alpha|\le k}(t^{-1}R)^{-\frac{k-|\alpha|}{d}}\biggr]\notag\\
 & & \qquad\qquad\qquad\qquad\qquad\qquad\qquad\qquad
 \times\int_{\{|y|\ge R^{1/d}\}}|y|^k|\varphi(y)|dy
 \label{eq:3.14}
 \end{eqnarray}
for all $t>0$ and $R>1$. 

Put $R=t+1$. 
Since 
$$
e^{t{\mathcal L}}\varphi(x)
=I_1(x,t)-I_2(x,t)+e^{t{\mathcal L}}\tilde{\psi}_R(x)+e^{t{\mathcal L}}\varphi_2(x), 
$$
by \eqref{eq:3.8}, \eqref{eq:3.9}, \eqref{eq:3.13} and \eqref{eq:3.14}
we have 
$$
 \int_{{\bf R}^N}|x|^\ell |e^{t{\mathcal L}}\varphi(x)|dx
 \preceq t^{-\frac{k-\ell}{d}}\int_{{\bf R}^N}|x|^k|\varphi(x)|dx
$$
for all $t>0$, and we obtain assertion~(ii). 
Furthermore, putting $R=\epsilon t+1$ with $\epsilon>0$, 
by \eqref{eq:3.8}, \eqref{eq:3.9}, \eqref{eq:3.13} and \eqref{eq:3.14} 
we have 
\begin{eqnarray*}
 & & \limsup_{t\to\infty}t^{\frac{k-\ell}{d}}\int_{{\bf R}^N}|x|^\ell |e^{t{\mathcal L}}\varphi(x)|dx
 \le\limsup_{t\to\infty}t^{\frac{k-\ell}{d}}\int_{{\bf R}^N} |x|^k|I_1(x,t)|dx\\
 & & \qquad\qquad
\le C_5\left[\epsilon^{\frac{[k]+1-k}{d}}+\epsilon^{\frac{[k]+1+\ell-k}{d}}\right]\int_{{\bf R}^N}|y|^k|\varphi(y)|dy
\end{eqnarray*}
for some constant $C_5$. 
Since $\epsilon$ is arbitrary, 
we obtain \eqref{eq:3.1}, and assertion~(iii) follows. 
Thus the proof of Proposition~\ref{Proposition:3.1} is complete.  
$\Box$
\vspace{5pt}
\newline
Now we are ready to prove Theorem~\ref{Theorem:1.1}.
\vspace{5pt}
\newline
{\bf Proof of Theorem~\ref{Theorem:1.1}.}
Let $v$ be the function given in Theorem~\ref{Theorem:1.1} and $0\le\ell\le K$. 
By Lemma~\ref{Lemma:2.2} we have 
\begin{equation}
\label{eq:3.15}
|||v(0)|||_K\le C|||\varphi|||_K
\end{equation}
for some constant $C$. 
Then, by Proposition~\ref{Proposition:3.1}~(i), \eqref{eq:2.3} and \eqref{eq:2.4} 
we have 
\begin{eqnarray}
 t^{\frac{N}{d}(1-\frac{1}{q})+\frac{j}{d}}\|\nabla^jv(t)\|_q
\!\!\! & \preceq &\!\!\!
\|v(t/2)\|_1\le \|v(0)\|_1\preceq|||\varphi|||_K,\label{eq:3.16}\\
 t^{\frac{j}{d}-\frac{\ell}{d}}|||\nabla^jv(t)|||_\ell
 \!\!\! & \preceq &\!\!\!
\|v(t/2)\|_1+t^{-\frac{\ell}{d}}\int_{{\bf R}^N}|x|^\ell |v(t/2)|dx,\label{eq:3.17}\\
t^{\frac{j}{d}}(1+t)^{-\frac{\ell}{d}}|||\nabla^jv(t)|||_\ell
\!\!\! & \preceq &\!\!\!
|||v(0)|||_K\preceq |||\varphi|||_K,\label{eq:3.18}
\end{eqnarray}
for all $t>0$. 
On the other hand, it follows from Lemma~\ref{Lemma:2.1} that 
$$
\int_{{\bf R}^N}x^\alpha v(x,0)dx=0,\quad \alpha\in{\bf M}_K.
$$
Therefore, applying Proposition~\ref{Proposition:3.1}~(ii) with the aid of \eqref{eq:3.15}, 
we see that 
\begin{eqnarray}
\|v(t)\|_1+t^{-\frac{\ell}{d}}\int_{{\bf R}^N}|x|^\ell |v(x,t)|dx
\preceq t^{-\frac{K}{d}}\int_{{\bf R}^N}|x|^K|v(x,0)|dx\preceq t^{-\frac{K}{d}}|||\varphi|||_K
\label{eq:3.19}
\end{eqnarray}
for all $t>0$. Similarly, by Proposition~\ref{Proposition:3.1}~(iii) 
we have  
\begin{equation}
\label{eq:3.20}
\lim_{t\to\infty}t^{\frac{K}{d}}
\biggr[\|v(t)\|_1+t^{-\frac{\ell}{d}}\int_{{\bf R}^N}|x|^\ell |v(x,t)|dx\biggr]=0. 
\end{equation}
Hence, by \eqref{eq:3.16}--\eqref{eq:3.20} we have 
\begin{eqnarray*}
 & & t^{\frac{N}{d}(1-\frac{1}{q})+\frac{j}{d}}\|\nabla^jv(t)\|_q+t^{\frac{j}{d}}(1+t)^{-\frac{\ell}{d}}|||\nabla^jv(t)|||_\ell\\
 & & \qquad\quad
 \le C_1\min\{|||\varphi|||_K,t^{-\frac{K}{d}}|||\varphi|||_K\}
 \le C_2 (1+t)^{-\frac{K}{d}}|||\varphi|||_K,\qquad t>0,
 \end{eqnarray*}
where $C_1$ and $C_2$ are constants independent of $\varphi\in L^1_K$, and 
$$
\lim_{t\to\infty}t^{\frac{K}{d}}
\biggr[t^{\frac{N}{d}(1-\frac{1}{q})+\frac{j}{d}}\|\nabla^jv(t)\|_q+t^{\frac{j}{d}}(1+t)^{-\frac{\ell}{d}}|||\nabla^jv(t)|||_\ell\biggr]=0.
$$
Thus we obtain \eqref{eq:1.12} and \eqref{eq:1.13}, 
and the proof of Theorem~\ref{Theorem:1.1} is complete. 
$\Box$
\section{Proof of Theorem~\ref{Theorem:1.2}}
In this section we prove Theorem~\ref{Theorem:1.2} 
by using Proposition~\ref{Proposition:3.1} and the operator $P_k(t)$, 
\vspace{5pt}
\newline
{\bf Proof of Theorem~\ref{Theorem:1.2}.}
Assertion~(i) follows from \eqref{eq:2.9}. 
We prove assertion~(ii). 
Let $1\le q\le\infty$ and $0\le \ell\le K$. 
For any $j\in\{0,\cdots,\gamma\}$ with $j<d$, 
put 
\begin{eqnarray*}
 & & I_1(t):=\int_{t/2}^t \nabla^je^{(t-s){\mathcal L}}P_K(s)f(s)ds,\\
 & & I_2(t):=\int_0^{t/2} \nabla^je^{(t-s){\mathcal L}}P_K(s)f(s)ds
 =\int_0^{t/2} \nabla^je^{\frac{t-s}{2}{\mathcal L}}e^{\frac{t-s}{2}{\mathcal L}}P_K(s)f(s)ds.
\end{eqnarray*}
By \eqref{eq:2.3} and \eqref{eq:2.10} we have 
\begin{eqnarray}
 & & t^{\frac{N}{d}(1-\frac{1}{q})}\|I_1(t)\|_q
\preceq t^{\frac{N}{d}(1-\frac{1}{q})}\int_{t/2}^t (t-s)^{-\frac{j}{d}}\|P_K(s)f(s)\|_q\,ds\notag\\
 & & \quad
 \preceq
 \int_{t/2}^t (t-s)^{-\frac{j}{d}}s^{\frac{N}{d}(1-\frac{1}{q})}\|P_K(s)f(s)\|_q\,ds
 \preceq \int_{t/2}^t(t-s)^{-\frac{j}{d}}(1+s)^{-\frac{K}{d}}E_{K,q}[f](s)ds\notag\\
  & & \quad
  \preceq (1+t)^{-\frac{K}{d}}\int_{t/2}^t(t-s)^{-\frac{j}{d}}E_{K,q}[f](s)ds
\label{eq:4.1}
\end{eqnarray}
for all $t>0$. 
Furthermore, applying Proposition~\ref{Proposition:3.1}~(i) with the aid of \eqref{eq:2.10}, 
we obtain
\begin{eqnarray}
 & & (1+t)^{-\frac{\ell}{d}}|||I_1(t)|||_\ell\notag\\
 & & \preceq (1+t)^{-\frac{\ell}{d}}\int_{t/2}^t\biggr[(t-s)^{-\frac{j-\ell}{d}}\|P_K(s)f(s)\|_1+(t-s)^{-\frac{j}{d}}|||P_K(s)f(s)|||_\ell\biggr]\,ds\notag\\
 & & 
\preceq\int_{t/2}^t(t-s)^{-\frac{j}{d}}\biggr[\|P_K(s)f(s)\|_1+(1+s)^{-\frac{\ell}{d}}|||P_K(s)f(s)|||_\ell\biggr]\,ds\notag\\
 & & 
 \preceq (1+t)^{-\frac{K}{d}}\int_{t/2}^t(t-s)^{-\frac{j}{d}}E_{K,q}[f](s)ds
\label{eq:4.2}
\end{eqnarray}
for all $t>0$. 

On the other hand, applying Proposition~\ref{Proposition:3.1}~(ii) with the aid of Lemma~\ref{Lemma:2.1}~(ii), 
for any $\delta>0$, 
we deduce from \eqref{eq:2.10} that 
\begin{equation}
\label{eq:4.3}
|||e^{\frac{t-s}{2}{\mathcal L}}P_K(s)f(s)|||_\ell
\preceq (t-s)^{-\frac{K-\ell}{d}}|||P_K(s)f(s)|||_K
\preceq (t-s)^{-\frac{K-\ell}{d}}E_{K,q}[f](s)
\end{equation}
for all $t\ge s+\delta>0$. 
Similarly to \eqref{eq:4.1} and \eqref{eq:4.2}, 
we have 
\begin{eqnarray}
 & & t^{\frac{N}{d}(1-\frac{1}{q})}\|I_2(t)\|_q+(1+t)^{-\frac{\ell}{d}}|||I_2(t)|||_\ell\notag\\
 & & \preceq t^{\frac{N}{d}(1-\frac{1}{q})}
\int_0^{t/2} (t-s)^{-\frac{N}{d}(1-\frac{1}{q})-\frac{j}{d}}\|e^{\frac{t-s}{2}{\mathcal L}}P_K(s)f(s)\|_1\,ds\notag\\
 & & \qquad\quad
 +(1+t)^{-\frac{\ell}{d}}\int_0^{t/2}(t-s)^{-\frac{j-\ell}{d}}\|e^{\frac{t-s}{2}{\mathcal L}}P_K(s)f(s)\|_1ds\notag\\
 & & \qquad\qquad\qquad
 +(1+t)^{-\frac{\ell}{d}}\int_0^{t/2}(t-s)^{-\frac{j}{d}}|||e^{\frac{t-s}{2}{\mathcal L}}P_K(s)f(s)|||_\ell\,ds\notag\\
 & & \preceq\int_0^{t/2}(t-s)^{-\frac{j}{d}}
 \biggr[\|e^{\frac{t-s}{2}{\mathcal L}}P_K(s)f(s)\|_1+(t-s)^{-\frac{\ell}{d}}|||e^{\frac{t-s}{2}{\mathcal L}}P_K(s)f(s)|||_\ell\biggr]\,ds
 \label{eq:4.4}
 \end{eqnarray}
 for all $t>0$. 
On the other hand, for any $T>0$, 
it follows from Proposition~\ref{Proposition:3.1}~(iii) and Lemma~\ref{Lemma:2.1}~(ii) that 
$$
\lim_{t\to\infty}(t-s)^{\frac{K-k}{d}}|||e^{\frac{t-s}{2}{\mathcal L}}P_K(s)f(s)|||_k=0
$$
for any $0\le k\le K$ and $s\in(0,T)$.  
Then, by the Lebesgue dominated convergence theorem and Proposition~\ref{Proposition:3.1}~(ii)
we see that 
\begin{eqnarray}
 & & \limsup_{t\to\infty}
t^{\frac{K+j}{d}}\int_0^T(t-s)^{-\frac{j}{d}}
 \biggr[\|e^{\frac{t-s}{2}{\mathcal L}}P_K(s)f(s)\|_1+(t-s)^{-\frac{\ell}{d}}|||e^{\frac{t-s}{2}{\mathcal L}}P_K(s)f(s)|||_\ell\biggr]\,ds
 \notag\\
 & & \le\limsup_{t\to\infty}
 \int_0^T(t-s)^{\frac{K}{d}}
 \biggr[\|e^{\frac{t-s}{2}{\mathcal L}}P_K(s)f(s)\|_1+(t-s)^{-\frac{\ell}{d}}|||e^{\frac{t-s}{2}{\mathcal L}}P_K(s)f(s)|||_\ell\biggr]\,ds
 \notag\\
 & & =0. 
\end{eqnarray}
Furthermore, by \eqref{eq:4.3}, for any $T_0>0$, we can find constant $C_1$ and $C_2$ such that  
\begin{eqnarray}
 & & \int_T^{t/2}(t-s)^{-\frac{j}{d}}
 \biggr[\|e^{\frac{t-s}{2}{\mathcal L}}P_K(s)f(s)\|_1+(t-s)^{-\frac{\ell}{d}}|||e^{\frac{t-s}{2}{\mathcal L}}P_K(s)f(s)|||_\ell\biggr]\,ds
 \notag\\
 & & \le C_1\int_T^{t/2}(t-s)^{-\frac{K}{d}-\frac{j}{d}}E_{K,q}[f](s)ds
 \le C_2 t^{-\frac{K}{d}}\int_T^{t/2}(t-s)^{-\frac{j}{d}}E_{K,q}[f](s)ds
\label{eq:4.6}
\end{eqnarray}
for all $t\ge 2T$ and $T\ge T_0$. 
Therefore, by \eqref{eq:4.4}--\eqref{eq:4.6}, 
for any $\epsilon>0$ and $T\ge T_0$, 
we have 
\begin{equation}
t^{\frac{N}{d}(1-\frac{1}{q})}\|I_2(t)\|_q+t^{-\frac{\ell}{d}}|||I_2(t)|||_\ell
\le\epsilon t^{-\frac{K+j}{d}}+C_3t^{-\frac{K}{d}}\int_T^{t/2}(t-s)^{-\frac{j}{d}}E_{K,q}[f](s)ds
\label{eq:4.7}
\end{equation}
for all sufficiently large $t$, where $C_3$ is a constant independent of $T\in[T_0,\infty)$ and $\epsilon>0$. 
Hence, 
by \eqref{eq:4.1}, \eqref{eq:4.2} and \eqref{eq:4.7} 
we have \eqref{eq:1.16}. 
In addition, \eqref{eq:1.17} immediately follows from \eqref{eq:1.16}.
Thus the proof of Theorem~\ref{Theorem:1.2} is complete. 
$\Box$
\section{Integral equation with power nonlinearity}
Let $F=F(x,t,u)$ be a function in ${\bf R}^N\times(0,\infty)\times{\bf R}$ 
such that 
\begin{eqnarray}
\label{eq:5.0}
 & & F(x,t,0)=0,\vspace{5pt}\\ 
\label{eq:5.1}
 & & 
 |F(x,t,u_1)-F(x,t,u_2)|\le C_*(1+t)^{A}\max\{|u_1|^{p-1},|u_2|^{p-1}\}|u_1-u_2|,
\end{eqnarray} 
for $x\in{\bf R}^N$, $t>0$ and $u_1$, $u_2\in{\bf R}$, where $C_*>0$, $A\in{\bf R}$ and $p\ge 1$. 
Consider the integral equation
\begin{equation}
\label{eq:5.2}
u(x,t)=\int_{{\bf R}^N}G(x-y,t)\varphi(y)dy+\int_0^t\int_{{\bf R}^N}G(x-y,t-s)F(y,s,u(y,s))dyds,
\end{equation} 
where $\varphi\in L^1_K$ for some $K\ge 0$. 
Problem~\eqref{eq:5.2} is a generalization of problems~\eqref{eq:1.4} and \eqref{eq:1.5}. 
In this section, under condition~$({\bf G})$ and \eqref{eq:5.1}, 
we study the asymptotic behavior of the solution $u$ of \eqref{eq:5.2} satisfying 
\begin{equation}
\label{eq:5.3}
\sup_{t>0}\,(1+t)^{\frac{N}{d}(1-\frac{1}{q})}\|u(t)\|_q+\sup_{t>0}\,(1+t)^{-\frac{\ell}{d}}|||u(t)|||_\ell<\infty
\end{equation}
for any $q\in[1,\infty]$ and $\ell\in[0,K]$, 
and prove the following theorem. 
(For the existence of the solutions of \eqref{eq:5.2} satisfying \eqref{eq:5.3}, 
see \cite{IKK01}.) 
\begin{theorem}
\label{Theorem:5.1}
Assume condition $({\bf G})$ for some $\gamma\in{\bf N}$, $d>0$ and $L>0$. 
Let $0\le K<L$ with $[K]+1\le\gamma$ and $\varphi\in L^\infty\cap L^1_K$. 
Assume \eqref{eq:5.0}, \eqref{eq:5.1} and 
$$
A_p:=-A+N(p-1)/d-1>0.
$$ 
Let $u$ be a global-in-time solution of \eqref{eq:5.2} satisfying \eqref{eq:5.3}. 
\vspace{3pt}
\newline
{\rm (i)} For any $\alpha\in{\bf M}_K$, put 
\begin{equation}
\label{eq:5.4}
c_\alpha(t):=M_\alpha(\varphi,0)+\int_0^t M_\alpha(F(s),s)ds,
\end{equation}
where $F(x,t):=F(x,t,u(x,t))$. 
If  $A_p>|\alpha|/d$, then there exists a constant $c_\alpha$ such that 
\begin{equation}
\label{eq:5.5-1}
|c_\alpha(t)-c_\alpha|=O(t^{-A_p+|\alpha|/d})
\end{equation}
as $t\to\infty$. 
If $A_p\le|\alpha|/d$, then 
\begin{equation}
\label{eq:5.5}
c_\alpha(t)=\left\{
\begin{array}{ll}
O(t^{-A_p+|\alpha|/d}) & \mbox{if}\quad A_p<|\alpha|/d,\vspace{3pt}\\
O(\log t) & \mbox{if}\quad A_p=|\alpha|/d,
\end{array}
\right.
\end{equation}
as $t\to\infty$. 
\vspace{3pt}
\newline
{\rm (ii)} 
Define the functions $U_n=U_n(x,t)$ $(n=0,1,2,\dots)$ inductively by 
\begin{eqnarray}
U_0(x,t)\!\!\! & := &\!\!\!\sum_{|\alpha|\le K}c_\alpha(t)g_\alpha(x,t),\label{eq:5.6}\\
U_n(x,t) \!\!\! & := &\!\!\!
U_0(x,t)+\int_0^t e^{(t-s){\mathcal L}}P_K(s)F_{n-1}(s)ds\nonumber\\
 \!\!\! & \,\,= &\!\!\!
 \sum_{|\alpha|\le K}\left[M_\alpha (\varphi, 0)+\int_{0}^{t}
M_\alpha (F(s)-F_{n-1}(s),s)ds\right]g_\alpha(x,t)\notag\\
 & & \qquad\qquad\qquad\qquad\qquad\qquad\quad\,\,\,\,
+\int_0^t e^{(t-s){\mathcal L}}F_{n-1}(s)ds, 
\label{eq:5.7}
\end{eqnarray}
where $n=1,2,\dots$ and $F_{n-1}(x,t):=F(x,t,U_{n-1}(x,t))$. 
Then, for any $q\in[1,\infty]$ and $\ell\in[0,K]$, 
\begin{equation}
\sup_{t>0}\,t^{\frac{N}{d}(1-\frac{1}{q})}\|U_n(t)\|_q+\sup_{t>0}\,(1+t)^{-\frac{\ell}{d}}|||U_n(t)|||_\ell<\infty
\label{eq:5.8}
\end{equation}
and 
\begin{eqnarray}
 & & t^{\frac{N}{d}(1-\frac{1}{q})}\left\|u(t)-U_n(t)\right\|_q
 +t^{-\frac{\ell}{d}}
\left|\left|\left|u(t)-U_n(t)\right|\right|\right|_l\vspace{3pt}\notag\\
 & & \qquad\qquad\qquad
 =\left\{
\begin{array}{ll}
o(t^{-\frac{K}{d}})+O(t^{-(n+1)A_p})  & \mbox{if}\quad (n+1)A_p\not=K/d,\vspace{3pt}\\
O(t^{-\frac{K}{d}}\log t)  & \mbox{if}\quad (n+1)A_p=K/d,\\
\end{array}
\right.
\label{eq:5.9}
\end{eqnarray}
as $t\to\infty$. 
\end{theorem}
Here we remark: 
\begin{itemize}
  \item $U_0(\cdot,t)$ is a linear combination of $\{g_\alpha(\cdot,t)\}_{|\alpha|\le K}$ 
	and plays a role of the projection of $u(\cdot,t)$ into the finite dimensional space spanned 
	by $\{g_\alpha(\cdot,t)\}_{|\alpha|\le K}$;  
  \item For $n=1,2,\dots$, $U_n$ is a nonlinear approximation to the solution $u$ and 
	is constructed by $U_0$ systematically. 
\end{itemize} 
{\bf Proof of Theorem~\ref{Theorem:5.1}.} 
Let $f(x,t)=F(x,t,u(x,t))$. 
It follows from \eqref{eq:5.0} and \eqref{eq:5.1} that 
$$
|f(x,t)|\le C_*(1+t)^{A}|u(x,t)|^p. 
$$
This together with \eqref{eq:1.14} and \eqref{eq:5.3} implies 
\begin{eqnarray}
E_{K,q}[f](t)
\!\!\! & \le &\!\!\! C_*(1+t)^{A}\|u(t)\|_\infty^{p-1}E_{K,q}[u](t)\notag\\
\!\!\! & \preceq &\!\!\!
(1+t)^{A-\frac{N}{d}(p-1)+\frac{K}{d}}
=(1+t)^{-A_p-1+\frac{K}{d}}
\label{eq:5.10}
\end{eqnarray}
for all $t>0$. Then, by Lemma~\ref{Lemma:2.2} and \eqref{eq:5.4} we have 
$$
|c_\alpha(t_2)-c_\alpha(t_1)|
\le\int_{t_1}^{t_2}|M_\alpha(f(s),s)|ds
\preceq\int_{t_1}^{t_2} (1+s)^{-A_p-1+\frac{|\alpha|}{d}}ds
$$
for $t_2\ge t_1>0$. This implies 
\eqref{eq:5.5-1} and \eqref{eq:5.5}, and assertion~(i) follows. 

We prove assertion (ii). The proof is by induction. 
Assertion~(i) together with \eqref{eq:2.2} yields \eqref{eq:5.8} for $n=0$. 
Let $v=v(x,t)$ and $R_K[f](x,t)$ be functions given in Theorems~\ref{Theorem:1.1} and \ref{Theorem:1.2}. 
Since
\begin{eqnarray*}
 & & u(x,t)-U_0(x,t)\\
 & & =e^{t{\mathcal L}}\varphi(x)+\int_0^t e^{(t-s){\mathcal L}}f(s)ds
-\sum_{|\alpha|\le K}\biggr[M_\alpha(\varphi,0)+\int_0^t M_\alpha(f(s),s)ds\biggr]g_\alpha(x,t)\\
 & & =e^{t{\mathcal L}}\biggr[\varphi-\sum_{|\alpha|\le k}M_\alpha(\varphi,0)g_\alpha(0)\biggr]
 +\int_0^t e^{(t-s){\mathcal L}}\biggr[f(s)-\sum_{|\alpha|\le k}M_\alpha(f(s),s)g_\alpha(s)\biggr]ds\\
 & & =v(x,t)+R_K[f](x,t), 
\end{eqnarray*}
applying Theorems~\ref{Theorem:1.1} and \ref{Theorem:1.2} with the aid of \eqref{eq:5.10}, 
we obtain
$$
t^{\frac{N}{d}(1-\frac{1}{q})}\|u(t)-U_0(t)\|_q
+t^{-\frac{\ell}{d}}|||u(t)-U_0(t)|||_\ell
=o(t^{-\frac{K}{d}})+O\biggr(t^{-\frac{K}{d}}\int_T^t s^{-A_p-1+\frac{K}{d}}ds\biggr)
$$
as $t\to\infty$, for any $T>0$. 
This together with  \eqref{eq:5.3} implies \eqref{eq:5.8} and \eqref{eq:5.9} for $n=0$.  
Thus assertion~(ii) holds for $n=0$. 

Assume that assertion~(ii) holds for some $n=m\in{\bf N}\cup\{0\}$. 
Then 
\begin{equation}
\label{eq:5.11}
u(x,t)-U_{m+1}(x,t)
=v(x,t)+\int_0^t e^{(t-s){\mathcal L}}P_K(s)f_m(s)ds,
\end{equation}
where $f_m(x,t):=F(x,t,u(x,t))-F(x,t,U_m(x,t))$. 
Similarly to \eqref{eq:5.10}, 
by \eqref{eq:5.1} and assertion~(ii) with $n=m$ 
we have 
\begin{eqnarray}
E_{K,q}[f_m](t)
\!\!\! & \preceq & \!\!\!
(1+t)^{A}\max\{\|u(t)\|_\infty^{p-1},\|U_m(t)\|_\infty^{p-1}\}E_{K,q}[u-U_m](t)\notag\\
\!\!\! & = &\!\!\! 
\left\{
\begin{array}{ll}
o(t^{-A_p-1})+O(t^{\frac{K}{d}-(m+2)A_p-1}) & \mbox{if}\quad (m+1)A_p\not=K/d,\vspace{3pt}\\
O(t^{-A_p-1}\log t) & \mbox{if}\quad (m+1)A_p=K/d,
\end{array}
\right.
\label{eq:5.12}
\end{eqnarray}
as $t\to\infty$. 
Then, by Theorem~\ref{Theorem:1.2}, for any $T>0$, 
we have
\begin{eqnarray}
 & & t^{\frac{N}{d}(1-\frac{1}{q})}\left\|\int_0^te^{(t-s){\mathcal L}}P_K(s)f_m(s)ds\right\|_q
+t^{-\frac{\ell}{d}}\biggr|\biggr|\biggr|\int_0^te^{(t-s){\mathcal L}}P_K(s)f_m(s)ds\biggr|\biggr|\biggr|_\ell\notag\\
 & & \qquad\qquad
 =o(t^{-\frac{K}{d}})+O\biggr(t^{-\frac{K}{d}}\int_T^t E_{K,q}[f_m](s)ds\biggr)\notag\\
 & & \qquad\qquad
 =\left\{
\begin{array}{ll}
o(t^{-\frac{K}{d}})+O(t^{-(m+2)A_p}) & \mbox{if}\quad (m+2)A_p\not=K/d,\vspace{3pt}\\
O(t^{-\frac{K}{d}}\log t) & \mbox{if}\quad (m+2)A_p=K/d,
\end{array}
\right.
\label{eq:5.13}
\end{eqnarray}
as $t\to\infty$. Therefore we deduce from Theorem~\ref{Theorem:1.1}, \eqref{eq:5.11} and \eqref{eq:5.13} that 
\begin{eqnarray*}
 & & t^{\frac{N}{d}(1-\frac{1}{q})}\|u(t)-U_{m+1}(t)\|_q
+t^{-\frac{\ell}{d}}|||u(t)-U_{m+1}(t)|||_\ell\\
 & & \qquad\qquad
=\left\{
\begin{array}{ll}
o(t^{-\frac{K}{d}})+O(t^{-(m+2)A_p}) & \mbox{if}\quad (m+2)A_p\not=K/d,\vspace{3pt}\\
O(t^{-\frac{K}{d}}\log t) & \mbox{if}\quad (m+2)A_p=K/d,
\end{array}
\right.
\end{eqnarray*}
as $t\to\infty$.
This together with \eqref{eq:5.3} implies \eqref{eq:5.8} and \eqref{eq:5.9} with $n=m+1$. 
Hence, by induction we see that assertion (ii) holds for all $n=0,1,2,\dots$, 
and the proof of Theorem~\ref{Theorem:5.1} is complete.
$\Box$
\vspace{5pt}
\newline
As a corollary of Theorem~\ref{Theorem:5.1} with $n=0$, 
we give a decay estimate of the distance in $L^q$ $(1\le q\le\infty)$ and $L^1_\ell$ $(0\le\ell\le K)$ 
from the solution $u$ of \eqref{eq:5.2} to $Mg(x,t)$, 
where 
$$
M:=\lim_{t\to\infty}\int_{{\bf R}^N}u(x,t)dx
=\int_{{\bf R}^N}\varphi(x)dx+\int_0^\infty\int_{{\bf R}^N}F(x,t,u(x,t))dxdt. 
$$
We remark that $M$ coincides with $c_0$, which is given in Theorem~\ref{Theorem:5.1}~(i). 
\begin{corollary}
\label{Corollary:5.1}
Assume the same conditions as in Theorem~{\rm\ref{Theorem:5.1}}. 
\newline
{\rm (i)} For any $1\le q\le\infty$ and $0\le\ell\le K$, 
\begin{eqnarray}
 & & t^{\frac{N}{d}(1-\frac{1}{q})}
\|u(t)-Mg(t)\|_q+t^{-\frac{\ell}{d}}|||u(t)-Mg(t)|||_\ell\vspace{5pt}\notag\\
 & & \qquad
=
\left\{
 \begin{array}{ll}
 o(t^{-\frac{K}{d}})+O(t^{-A_p}) & \mbox{if}\quad\mbox{$0\le K<1$ and $A_p\not=K/d$},\vspace{3pt}\\
 O(t^{-\frac{K}{d}}\log t) & \mbox{if}\quad\mbox{$0\le K<1$ and $A_p=K/d$},\vspace{3pt}\\
 O(t^{-\frac{1}{d}})+O(t^{-A_p}) & \mbox{if}\quad\mbox{$K\ge 1$ and $A_p\not=1/d$},\vspace{3pt}\\
 O(t^{-\frac{1}{d}}\log t) & \mbox{if}\quad\mbox{$K\ge 1$ and $A_p=1/d$},
 \end{array}
 \right.
\label{eq:5.14}
\end{eqnarray}
as $t\to\infty$.
\vspace{3pt}
\newline
{\rm (ii)} Let $K\ge 1$. Let $f_M(x,t)=F(x,t,Mg(x,t))$, and assume that 
\begin{equation}
\label{eq:5.15}
\int_0^\infty\left|M_\alpha(f_M(t),t)\right|dt<\infty
\end{equation}
for any $\alpha\in{\bf M}$ with $|\alpha|=1$.
Then, for any $1\le q\le\infty$ and $0\le\ell\le K$,
\begin{equation}
\label{eq:5.16}
t^{\frac{N}{d}(1-\frac{1}{q})}
\|u(t)-Mg(t)\|_q+t^{-\frac{\ell}{d}}|||u(t)-Mg(t)|||_\ell
=O(t^{-\frac{1}{d}})+O(t^{-A_p})
\end{equation}
as $t\to\infty$. 
\end{corollary}
{\bf Proof.} 
It follows from \eqref{eq:5.6} that 
$$
U_0(x,t)-Mg(x,t)
=(c_0(t)-M)g(x,t)+\sum_{1\le|\alpha|\le K}c_\alpha(t)g_\alpha(x,t).
$$
This together with Theorem~\ref{Theorem:5.1}~(i) and \eqref{eq:2.2} implies that 
\begin{eqnarray}
 & & t^{\frac{N}{d}(1-\frac{1}{q})}\|U_0(t)-Mg(t)\|_q+t^{-\frac{\ell}{d}}|||U_0(t)-Mg(t)|||_\ell\notag\\
 & & 
 \preceq |c_0(t)-M|
 +\sum_{1\le|\alpha|\le K}(1+t)^{-\frac{|\alpha|}{d}}|c_\alpha(t)|\notag\\
 & & 
 =\left\{
 \begin{array}{ll}
 O(t^{-A_p}) & \mbox{if}\quad 0\le K<1,\\
 O(t^{-\frac{1}{d}})+O(t^{-A_p}) & \mbox{if}\quad\mbox{$K\ge 1$ and $A_p\not=1/d$},\\
 O(t^{-\frac{1}{d}}\log t) & \mbox{if}\quad\mbox{$K\ge 1$ and $A_p=1/d$},
 \end{array}
 \right.\label{eq:5.17}
\end{eqnarray}
as $t\to\infty$. 
Combining \eqref{eq:5.17} with Theorem~\ref{Theorem:5.1}~(ii), we see that 
\begin{eqnarray*}
 & & t^{\frac{N}{d}(1-\frac{1}{q})}
\|u(t)-Mg(t)\|_q+t^{-\frac{\ell}{d}}|||u(t)-Mg(t)|||_\ell\vspace{5pt}\notag\\
 & & \qquad
=
\left\{
 \begin{array}{ll}
 o(t^{-\frac{K}{d}})+O(t^{-A_p}) & \mbox{if}\quad\mbox{$0\le K<1$ and $A_p\not=K/d$},\\
 O(t^{-\frac{K}{d}}\log t) & \mbox{if}\quad\mbox{$0\le K<1$ and $A_p=K/d$},\\
 O(t^{-\frac{1}{d}})+O(t^{-A_p}) & \mbox{if}\quad\mbox{$K\ge 1$ and $A_p\not=1/d$},\\
 O(t^{-\frac{1}{d}}\log t) & \mbox{if}\quad\mbox{$K\ge 1$ and $A_p=1/d$},\\
 \end{array}
 \right.
\end{eqnarray*}
as $t\to\infty$, and assertion~(i) follows. 

We prove assertion~(ii). It suffices to consider the case $A_p=1/d$. 
Let $K\ge 1$ and $\alpha\in{\bf M}$ with $|\alpha|=1$. 
By \eqref{eq:5.4} and \eqref{eq:5.15} 
we apply Lemmas~\ref{Lemma:2.1} and \ref{Lemma:2.2} to obtain 
\begin{equation*}
\begin{split}
|c_\alpha(t_2)-c_\alpha(t_1)|
 & \preceq\,\biggr|\int_{t_1}^{t_2}\left[M_\alpha(f(s),s)-M_\alpha(f_M(s),s)\right]ds\biggr|+1\\
 & =\,\biggr|\int_{t_1}^{t_2}M_\alpha(f(s)-f_M(s),s)ds\biggr|+1\\
 & \preceq\int_{t_1}^{t_2}(1+s)^{-\frac{K}{d}+\frac{1}{d}}E_{K,q}[f-f_M](s)\,ds+1
\end{split}
\end{equation*}
for all $0<t_1<t_2$, where $f(x,t)=F(x,t,u(x,t))$. 
Then, by a similar argument to \eqref{eq:5.10} with the aid of \eqref{eq:5.14} we see that
$$
|c_\alpha(t_2)-c_\alpha(t_1)|\preceq\int_{t_1}^{t_2} s^{\frac{1}{d}-A_p-1}E_{K,q}[u-Mg](s)ds+1
\preceq\int_{t_1}^{t_2} s^{-A_p-1}\log s\,ds+1
$$
for all sufficiently large $t_1$ and $t_2$ with $t_1<t_2$. This implies that 
$|c_\alpha(t)|=O(1)$ as $t\to\infty$. 
Therefore, by the same argument as in the proof of assertion~(i) 
we have \eqref{eq:5.16}. 
Thus assertion~(ii) follows, and the proof of Corollary~\ref{Corollary:5.1} is complete. 
$\Box$ 
\vspace{5pt}
\newline
Next, applying Theorem~\ref{Theorem:5.1} with $n=1$, 
we give more precise description of the asymptotic behavior of the solution of \eqref{eq:5.2} than in Corollary~\ref{Corollary:5.1}.  
\begin{corollary}
\label{Corollary:5.2}
Assume the same conditions as in Theorem~{\rm\ref{Theorem:5.1}} 
and $0\le K<1$. 
Let 
$$
f(x,t):=F(x,t,u(x,t)),\qquad
f_M(x,t):=F(x,t,Mg(x,t)),
$$
and put 
$$
\tilde{u}(x,t):=M'g(x,t)+\int_0^t e^{(t-s){\mathcal L}}f_M(s)ds,
$$
where
$$
M':=\int_{{\bf R}^N}\varphi(x)dx+\int_0^\infty\int_{{\bf R}^N}[f(x,t)-f_M(x,t)]dxdt.
$$
Then, for any $q\in[1,\infty]$ and $\ell\in[0,K]$, 
\begin{eqnarray}
 & & t^{\frac{N}{d}(1-\frac{1}{q})}
\|u(t)-\tilde{u}(t)\|_q+t^{-\frac{\ell}{d}}|||u(t)-\tilde{u}(t)|||_\ell\vspace{5pt}\notag\\
 & & \qquad
=
\left\{
\begin{array}{ll}
o(t^{-\frac{K}{d}})+O(t^{-2A_p}) & \mbox{if}\quad 2A_p\not=K/d,\vspace{3pt}\\
O(t^{-\frac{K}{d}}\log t) & \mbox{if}\quad 2A_p=K/d,\\
\end{array}
\right.
\label{eq:5.18}
\end{eqnarray}
as $t\to\infty$.  
\end{corollary}
{\bf Proof.}
Put 
\begin{eqnarray*}
& & f_1(x,t):=F(U_0(x,t))-f_M(x,t),\quad
f_2(x,t):=f(x,t)-f_M(x,t),\\
 & & w(x,t):=\int_0^t e^{(t-s){\mathcal L}}P_K(s)f_1(s)ds. 
\end{eqnarray*}
Similarly to \eqref{eq:5.12}, 
by Corollary~\ref{Corollary:5.1} and \eqref{eq:5.17} we have 
\begin{equation}
E_{K,q}[f_1](t)+E_{K,q}[f_2](t)
=\left\{
 \begin{array}{ll}
 o(t^{-A_p-1})+O(t^{\frac{K}{d}-2A_p-1}) & \mbox{if}\quad A_p\not=K/d,\vspace{3pt}\\
 O(t^{-A_p-1}\log t) & \mbox{if}\quad A_p=K/d,
 \end{array}
 \right.
 \label{eq:5.19}
\end{equation}
as $t\to\infty$.
Then, by Theorem~\ref{Theorem:1.2} and \eqref{eq:5.19}
we see that 
\begin{equation}
t^{\frac{N}{d}(1-\frac{1}{q})}\|w(t)\|_q
+t^{-\frac{\ell}{d}}|||w(t)|||_\ell
 =\left\{
 \begin{array}{ll}
 o(t^{-\frac{K}{d}})+O(t^{-2A_p}) & \mbox{if}\quad 2A_p\not=K/d,\vspace{3pt}\\
 O(t^{-\frac{K}{d}}\log t) & \mbox{if}\quad 2A_p=K/d,
 \end{array}
 \right.
 \label{eq:5.20}
\end{equation}
as $t\to\infty$. 
Furthermore, by Lemma~\ref{Lemma:2.2} and \eqref{eq:5.19} we get
\begin{eqnarray}
 & & |M_0(f_2(t),t)|\preceq
 (1+t)^{-\frac{K}{d}}E_{K,q}[f_2](t)\notag\\
 & & \qquad\qquad\qquad
 =
\left\{
 \begin{array}{ll}
 o(t^{-\frac{K}{d}-A_p-1})+O(t^{-2A_p-1}) & \mbox{if}\quad A_p\not=K/d,\vspace{3pt}\\
 O(t^{-2A_p-1}\log t) & \mbox{if}\quad A_p=K/d,
 \end{array}
 \right.\label{eq:5.21}
\end{eqnarray}
as $t\to\infty$. 
On the other hand, it follows from \eqref{eq:5.7} that 
$$
 \tilde{u}(x,t)-U_1(x,t)=
 \int_t^\infty M_0(f_2(s),s)ds\,\cdot g(x,t)-w(x,t).
$$
Therefore, by \eqref{eq:5.20} and \eqref{eq:5.21} we have 
\begin{eqnarray*}
 & & t^{\frac{N}{d}(1-\frac{1}{q})}\|U_1(t)-\tilde{u}(t)\|_q+t^{-\frac{\ell}{d}}|||U_1(t)-\tilde{u}(t)|||_\ell\\
 & & \qquad
 =\left\{
 \begin{array}{ll}
 o(t^{-\frac{K}{d}})+O(t^{-2A_p}) & \mbox{if}\quad 2A_p\not=K/d,\vspace{3pt}\\
 O(t^{-\frac{K}{d}}\log t) & \mbox{if}\quad 2A_p=K/d,
 \end{array}
 \right.
\end{eqnarray*}
as $t\to\infty$. This together with \eqref{eq:5.9} implies \eqref{eq:5.18}, 
and Corollary~\ref{Corollary:5.2} follows. 
$\Box$
\section{Applications}
We apply the results in the previous sections to some nonlinear parabolic equations, 
and show the validity of our arguments. 
\subsection{Semilinear parabolic equations}
Let $u$ be a solution of the Cauchy problem for a semilinear parabolic equation
\begin{equation}
\label{eq:6.1}
\left\{
\begin{array}{ll}
\partial_t u=\Delta u+a(x,t)|u|^{p-1}u & \mbox{in}\quad{\bf R}^N\times(0,\infty),\vspace{3pt}\\
u(x,0)=\varphi(x) & \mbox{in}\quad{\bf R}^N,
\end{array}
\right.
\end{equation}
where $p\ge 1$, $a\in L^\infty(0,\infty:L^\infty({\bf R}^N))$ and $\varphi\in L^\infty$. 
Asymptotic behavior of the solutions of \eqref{eq:6.1} has been studied by many mathematicians  
(see e.g. \cite{DK}, \cite{GV}, \cite{H}--\cite{IK03}, \cite{IKoba}, \cite{KP}, \cite{QS}, \cite{T} and references therein). 
In particular, the asymptotic expansions of the solutions of \eqref{eq:6.1} behaving like a multiple of the heat kernel 
were discussed in \cite{IIK}--\cite{IK03}. 

On the other hand, 
the heat kernel satisfies condition~$({\bf G})$ for any $\gamma\in{\bf N}$ and $L>0$ with $d=2$. 
Therefore, as a corollary of the results in the previous section, we have: 
\begin{theorem}
\label{Theorem:6.1}
Let $\varphi\in L^\infty\cap L^1_k$ for some $K\ge 0$ and $p\ge 1$. 
Assume 
\begin{equation}
\label{eq:6.2}
\sup_{t>0}\,(1+t)^{-A}\|a(t)\|_\infty<\infty
\end{equation}
for some $A\in{\bf R}$ and $A_p:=-A+N(p-1)/2-1>0$. 
Let $u$ be a solution of \eqref{eq:6.1} satisfying 
\begin{equation}
\label{eq:6.3}
\sup_{t>0}\,(1+t)^{\frac{N}{2}}\|u(t)\|_\infty<\infty.
\end{equation}
Then the conclusions as Theorem~{\rm\ref{Theorem:5.1}} 
and Corollaries~{\rm\ref{Corollary:5.1}} and {\rm\ref{Corollary:5.2}} hold 
for any $\gamma\in{\bf N}$ and $L>0$ with $d=2$. 
\end{theorem}
{\bf Proof.}
Under assumptions~\eqref{eq:6.2} and \eqref{eq:6.3}, 
by a similar argument as in \cite[Theorem~3.1]{IK03} we see that 
$$
\sup_{t>0}\,t^{\frac{N}{d}(1-\frac{1}{q})}\|u(t)\|_q+\sup_{t>0}\,(1+t)^{-\frac{\ell}{d}}|||u(t)|||_\ell<\infty
$$
for any $q\in[1,\infty]$ and $\ell\in[0,K]$. 
Then, applying the arguments in Section~5, 
we obtain the desired conclusions. Thus Theorem~\ref{Theorem:6.1} follows.
$\Box$
\begin{remark}
\label{Remark:6.1}
Let $G=G(x,t)$ be the heat kernel and $\alpha\in{\bf M}$ with $|\alpha|=1$. 
Then, by the radially symmetry  of $G$, 
we have 
$$
M_\alpha(f,t)=\int_{{\bf R}^N}x^\alpha f(x)dx,\qquad t>0,
$$
for all $f\in L^1_1$. 
Furthermore, if $a=a(x,t)$ is radially symmetric with respect to the space variable $x$, 
then 
$$
M_\alpha(f_M(t),t)=\int_{{\bf R}^N}x^\alpha f_M(x,t)dx
=|M|^{p-1}M\int_{{\bf R}^N}x^\alpha a(x,t)g(x,t)^pdx=0,
$$
where $f_M$ is the function defined in Corollary~{\rm\ref{Corollary:5.2}}, and 
assumption~\eqref{eq:5.15} is satisfied. 
\end{remark}
Theorem~\ref{Theorem:6.1} gives sharper decay estimates of 
$L^q({\bf R}^N)$-distance from the solution $u$ to its asymptotic profiles than 
in \cite{DK}, \cite{IIK}, \cite{IK01}--\cite{IK03} and \cite{T}. 
Furthermore, similarly to \cite{IK03}, 
we see that similar results to Theorem~\ref{Theorem:6.1} hold for 
more general nonlinear heat equations 
$$
\partial_t u=\Delta u+F(x,t,u,\nabla u)
\quad\mbox{in}\quad{\bf R}^N\times(0,\infty),
$$ 
under suitable assumptions on $F$ (see conditions~$(C_A)$ and $(F_A)$ in \cite{IK03}). 
The details are left to the reader. 
\subsection{Fractional semilinear parabolic equations}
Consider the Cauchy problem for a fractional semilinear parabolic equation 
\begin{equation}
\label{eq:6.4}
\left\{
\begin{array}{ll}
\partial_t u=-(-\Delta)^{\theta/2}u+a(x,t)|u|^{p-1}u & \mbox{in}\quad{\bf R}^N\times(0,\infty),\vspace{3pt}\\
u(x,0)=\varphi(x) & \mbox{in}\quad{\bf R}^N,
\end{array}
\right.
\end{equation}
where $0<\theta<2$, $p\ge 1$, $a\in L^\infty(0,\infty:L^\infty({\bf R}^N))$ and $\varphi\in L^\infty$.   
A continuous function $u$ in ${\bf R}^N\times(0,\infty)$ is said to be a solution of \eqref{eq:6.4} 
if $u$ satisfies 
$$
u(x,t)=\int_{{\bf R}^N}G_\theta(x-y,t)\varphi(y)dy+\int_0^t\int_{{\bf R}^N}G_\theta(x-y,t-s)a(y,s)|u(y,s)|^{p-1}u(y,s)dyds
$$
for all $x\in{\bf R}^N\times(0,\infty)$, where 
$G_\theta=G_\theta(x,t)$ be the fundamental solution of 
$$
\partial_t u+(-\Delta)^{\theta/2}u=0\qquad\mbox{in}\quad{\bf R}^N\times(0,\infty). 
$$
Problem~\eqref{eq:6.4} has been studied extensively by many mathematicians 
in view of various aspects, for example, 
nonlinear problems with anomalous diffusion and the Laplace equation with dynamical boundary conditions 
(see \cite{AF}, \cite{FIK}, \cite{FK}, \cite{HKN}, \cite{S}, \cite{Y02} and references therein). 
Among others, in the case where $a(x,t)$ is a negative constant function in ${\bf R}^N\times(0,\infty)$,
Fino and Karch \cite{FK} proved the following (see also \cite{HKN}):  
\begin{itemize}
  \item Let $\varphi\in L^1({\bf R}^N)$ and $p>1+\theta/N$. 
  Then there exists a constant  $M$ such that 
  \begin{equation}
  \label{eq:6.5}
  \lim_{t\to\infty}\int_{{\bf R}^N}u(x,t)dx=M
  \quad
  \mbox{and}
  \quad
  \lim_{t\to\infty}t^{\frac{N}{\theta}(1-\frac{1}{q})}\|u(t)-MG_\theta(t)\|_q=0
  \end{equation}
  for any $q\in[1,\infty]$. If $1<p\le 1+\theta/N$, then 
  $$
  \lim_{t\to\infty}\int_{{\bf R}^N}u(x,t)dx=0. 
  $$ 
\end{itemize}
(For the case $1<p\le 1+\theta/N$, see \cite{HKN}.) 
In the case where $a(x,t)$ is a positive constant function in ${\bf R}^N\times(0,\infty)$, 
the following holds: 
\begin{itemize}
  \item If $1<p\le 1+N/\theta$, then problem~\eqref{eq:6.4} has no positive global in time solutions (see~\cite{S}); 
  \item Let $\varphi\in L^\infty\cap L^1$, $\theta=1$ and $p>1+N$. 
  If $\|\varphi\|_1\|\varphi\|_\infty^{N(p-1)-1}$ is sufficiently small, 
  then there exists a global in time solution $u$ of \eqref{eq:6.4} such that \eqref{eq:6.5} holds with $\theta=1$ 
  (see~\cite{FIK}).
\end{itemize}
As far as we know, 
there are few results giving the precise description 
of the asymptotic behavior of the global in time solutions of \eqref{eq:6.4}. 
\vspace{3pt}

On the other hand, 
$G_\theta$ satisfies condition (${\bf G}$) for $d=L=\theta$ with $\gamma=1$ if $0<\theta\le 1$ and $\gamma=2$ if $1<\theta<2$ 
(see \cite[Lemma~7.3]{BFW},  \cite[Lemma~5.3]{BW} and \cite[Lemma~2.1]{Y02}). 
Then we can apply the results in Section~5 to problem~\eqref{eq:6.4}, 
and obtain the following theorem.
\begin{theorem}
\label{Theorem:6.2}
Consider problem~\eqref{eq:6.4}. Assume \eqref{eq:6.3} with $A_p:=-A+N(p-1)/\theta-1>0$. 
Then the conclusions of Theorem~{\rm\ref{Theorem:5.1}}  
and Corollaries~{\rm\ref{Corollary:5.1}} and {\rm\ref{Corollary:5.2}} hold for $d=L=\theta$ 
with $\gamma=1$ if $0<\theta\le 1$ and $\gamma=2$ if $1<\theta<2$. 
\end{theorem}
This enables us to study the precise description of the asymptotic behavior of the solution of \eqref{eq:6.4} 
behaving like a multiple of $G_\theta$ as $t\to\infty$ for the case $p>1+\theta/N$, 
and improves \cite{FIK} and \cite{FK}. 
\begin{remark}
\label{Remark:6.2}
Yamamoto~{\rm\cite{Y02}} recently studied the asymptotic behavior of the solutions of 
\begin{equation}
\label{eq:6.6}
\partial_t u=-(-\Delta)^{\theta/2}u+a(x,t)u\qquad
\mbox{in}\quad{\bf R}^N\times(0,\infty),
\qquad 1<\theta<2,
\end{equation}
and obtained higher order asymptotic expansions of the solutions, 
which are similar to those given in Theorem~{\rm\ref{Theorem:6.2}} with $p=1$ and $1<\theta<2$. 
However, his results require a stronger assumption on $a=a(x,t)$ than \eqref{eq:6.2} 
and a pointwise decay condition of the solution as $|x|\to\infty$, such as \eqref{eq:1.6}. 
\end{remark}
\subsection{Higher-order semilinear parabolic equations}
Let $m=1,2,\dots$ and 
$$
Lu:=\sum_{|\alpha|=2m}A_\alpha\partial_x^\alpha u
$$
be a $2m$-th order differential operator such that 
\begin{equation}
\label{eq:6.7}
\sum_{|\alpha|=2m}(i\xi)^\alpha A_\alpha\le -c_1|\mbox{Re}\,\xi|^{2m}+c_2|\mbox{Im}\,\xi|^m,\qquad\xi\in{\bf C}^N,
\end{equation}
for some positive constants $c_1$ and $c_2$, where $\{A_\alpha\}\subset{\bf R}$. 
In this section, 
under assumptions~\eqref{eq:5.1} and \eqref{eq:6.7}, 
we consider the Cauchy problem for the $2m$-th order semilinear parabolic equation 
\begin{equation}
\label{eq:6.8}
\left\{
\begin{array}{ll}
\partial_t u=Lu+a(x,t)|u|^p & 
\mbox{in}\quad{\bf R}^N\times(0,\infty),\vspace{3pt}\\
u(x,0)=\varphi(x) & \mbox{in}\quad{\bf R}^N,
\end{array}
\right.
\end{equation}
where $p\ge 1$, $a\in L^\infty(0,\infty:L^\infty({\bf R}^N))$ and $\varphi\in L^\infty\cap L^1$. 
In the case where $a$ is a positive constant function in ${\bf R}^N\times(0,\infty)$, 
problem~\eqref{eq:6.8} has been studied in several papers 
(see \cite{CM}, \cite{EGKP1}--\cite{EGKP3}, \cite{GP},  \cite{GG} and references therein), 
and the following holds:  
\begin{itemize}
 \item Let $1<p\le 1+2m/N$. If $\varphi\not\equiv 0$ in ${\bf R}^N$ and $\int_{{\bf R}^N}\varphi(x)dx\ge 0$, 
  then problem~\eqref{eq:6.8} has no global in time solutions (see \cite{EGKP1}); 
\item Let $p>1+2m/N$. Assume that $\varphi\not\equiv 0$ in ${\bf R}^N$ and $\int_{{\bf R}^N}\varphi(x)dx\ge 0$. 
 Then there exists a positive constant $C_1$ such that, if  
 $$
 |\varphi(x)|\le C_1e^{-|x|^{2m/(2m-1)}}
 \quad\mbox{for almost all $x\in{\bf R}^N$},
 $$
 then problem~\eqref{eq:6.8} has a global in time solution 
 behaving like a multiple of $G_m(x,t)$ as $t\to\infty$, where 
 $G_m=G_m(x,t)$ is the fundamental solution of 
 $$
 \partial_t u+(-\Delta)^m u=0\qquad\mbox{in}\quad{\bf R}^N\times(0,\infty)$$
 (see \cite{GP});
 \item Let $p>1+2m/N$. Assume 
 $$
 0\le\varphi(x)\le\frac{C_2}{1+|x|^\beta}
 $$
 for some $\beta>2m/(p-1)$ and $C_2>0$. If $\|\varphi\|_\infty$ is sufficiently small, 
 then problem~\eqref{eq:6.8} has a global in time solution 
 (see \cite{CM}). For the case $\beta=2m/(p-1)$, see \cite{GG}. 
\end{itemize}
Similarly to problem~\eqref{eq:6.4}, as far as we know, 
there are few results giving the precise description 
of the asymptotic behavior of the global in time solutions of  \eqref{eq:6.8}. 

On the other hand, under  assumption~\eqref{eq:6.7}, 
the fundamental solution of $\partial_t u=Lu$ in ${\bf R}^N\times(0,\infty)$ satisfies 
condition~$({\bf G})$ for any $\gamma\in{\bf N}$ and $L>0$ with $d=2m$ (see e.g.~\cite{Cui}). 
Then we apply  the results in Section~5 to problem~\eqref{eq:6.8}, 
and obtain the following theorem. 
\begin{theorem}
\label{Theorem:6.3}
Assume \eqref{eq:6.7}, 
and consider problem~\eqref{eq:6.8}. 
Assume \eqref{eq:6.3} with $A_p:=-A+N(p-1)/2m-1>0$. 
Then the conclusions of Theorem~{\rm\ref{Theorem:5.1}}  
and Corollaries~{\rm\ref{Corollary:5.1}} and {\rm\ref{Corollary:5.2}} hold for any $L>0$ and $\gamma>0$ with $d=2m$. 
\end{theorem}
Theorem~\ref{Theorem:6.3} enables us to study the precise description of the asymptotic behavior of the solutions 
behaving like a multiple of the kernel $G_m$. 
\bibliographystyle{amsplain}

\end{document}